\documentclass[11pt]{amsart}
\usepackage[margin=1.05in]{geometry}
\usepackage{amsmath,amssymb,amsthm,enumitem,graphicx,booktabs,xcolor,mathrsfs}
\usepackage{tikz}
\usetikzlibrary{arrows.meta,decorations.pathreplacing}
\definecolor{hl}{HTML}{EB6834}\definecolor{bl}{HTML}{2A78D6}\definecolor{gr}{HTML}{1BAF7A}
\usepackage[hypertexnames=false,colorlinks,linkcolor=blue!60!black,citecolor=green!40!black,urlcolor=blue!60!black]{hyperref}
\numberwithin{equation}{section}
\newtheorem{theorem}{Theorem}[section]
\newtheorem{proposition}[theorem]{Proposition}
\newtheorem{lemma}[theorem]{Lemma}

\newtheorem{maintheorem}{Theorem}
\newtheorem{maincorollary}[maintheorem]{Corollary}
\theoremstyle{definition}
\newtheorem{definition}[theorem]{Definition}
\theoremstyle{remark}
\newtheorem{remark}[theorem]{Remark}
\newtheorem{example}[theorem]{Example}
\newcommand{\F}{\mathbb F}\newcommand{\Z}{\mathbb Z}\newcommand{\Q}{\mathbb Q}\newcommand{\R}{\mathbb R}\newcommand{\C}{\mathbb C}
\newcommand{\eps}{\varepsilon}
\newcommand{\e}{\mathrm e}
\newcommand{\cE}{\mathscr E}

\title[Clusters of Farey fractions with polynomial denominators]{Clusters of Farey fractions with polynomial denominators\\ and the large sieve}

\author{Stephan Baier}
\address{Stephan Baier, Ramakrishna Mission Vivekananda Educational and Research Institute, Department of Mathematics, G.~T.~Road, PO Belur Math, Howrah, West Bengal 711202, India}
\email{stephanbaier2017@gmail.com, stephan@gm.rkmvu.ac.in}
\author{Chandler C. Corrigan}
\address{Chandler C. Corrigan, School of Mathematics and Statistics, University of New South Wales, Sydney NSW 2052, Australia}
\email{c.corrigan@unsw.edu.au}
\author{Sean B. Lynch}
\address{Sean B. Lynch, Westlake University, No.~600 Dunyu Road, Sandun Town, Xihu District, Hangzhou, Zhejiang Province 310030, China}
\email{sblynch@westlake.edu.cn}
\author{Liangyi Zhao}
\address{Liangyi Zhao, School of Mathematics and Statistics, University of New South Wales, Sydney NSW 2052, Australia}
\email{l.zhao@unsw.edu.au}

\date{Draft of 28 November 2026}

\begin{document}
	
	\begin{abstract}
		Let $f\in\Z[X]$ be a polynomial of degree $k\ge2$ with positive leading coefficient. For positive integers $Q$ and $N$ let $\Delta_f(Q,N)$ be the least number $\Delta$ such that the large sieve inequality
		\[
		\sum_{q\le Q}\ \sum_{\substack{a\bmod f(q)\\ \gcd(a,f(q))=1}}\Bigl|\sum_{M<n\le M+N}z_n\,\mathrm e^{2\pi\mathrm i\,an/f(q)}\Bigr|^2\le\Delta\sum_{M<n\le M+N}|z_n|^2
		\]
		holds for all integers $M$ and all complex numbers $z_n$. It has been expected that $\Delta_f(Q,N)\ll Q^\eps(Q^{k+1}+N)$. We show that this fails for every $f$ which is not a constant times a power of a linear polynomial. The reason is that the fractions $a/f(q)$ form clusters. If $\ell$ denotes the largest ramification index of $f$ over a non-zero complex value, we construct clusters which give $\Delta_f(Q,N)\gg Q^{1+k/\ell}N^{1-1/\ell}$; at $N=Q^{k+1}$ this exceeds the expected bound by the factor $Q^{1-1/\ell}\ge Q^{1/2}$. A second kind of cluster comes from approximating $f$ by its leading monomial. For $f=X^5+X+1$ and $N=Q^6$ it gives $\Delta_f(Q,N)\gg NQ^{3/5}$, whereas ramification gives $NQ^{1/2}$.
		
		In the opposite direction we prove $\Delta_f(Q,N)\ll Q^\eps(Q^{k+1}+N+Q^k\sqrt N)$ for every $f$. For quadratic polynomials of non-zero discriminant the lower and the upper bound agree up to the factor $Q^\eps$; this shows that a stronger estimate stated in the literature cannot hold. Finally, for the square moduli $q^2$ we prove $\Delta_f(Q,N)\ll Q^\eps(Q^3+N+Q^{8/5}N^{3/5})$. At $N=Q^3$ this is the first unconditional improvement of the bound of Baier and Zhao: the exponent of $Q$ by which the bound exceeds the expected one drops from $\frac12$ to $\frac25$. The new bound is better than that of Baier and Zhao for $Q^{11/4}<N<Q^4$.
	\end{abstract}
	
	\maketitle
	\tableofcontents
	
	\section{Introduction}\label{sec:intro}
	
	\subsection{The large sieve and sparse sets of moduli}
	Write $\e(x)=\exp(2\pi ix)$. The large sieve inequality, in the sharp form described in Montgomery's survey \cite{Mo}, says the following. If $x_1,\dots,x_R$ are real numbers whose pairwise distances modulo $1$ are at least $\delta>0$, then
	\begin{equation}\label{eq:LSgeneral}
		\sum_{i=1}^R\Bigl|\sum_{M<n\le M+N}z_n\,\e(nx_i)\Bigr|^2\le\bigl(\delta^{-1}+N-1\bigr)\sum_{M<n\le M+N}|z_n|^2
	\end{equation}
	for all integers $M$, all positive integers $N$ and all complex numbers $z_n$. The Farey fractions $a/s$ with $s\le Q$ have distances at least $Q^{-2}$. This gives the classical large sieve inequality, in which the factor on the right-hand side is $Q^2+N-1$.
	
	In many applications one needs such an inequality for a thin set of denominators, and the most natural thin sets are the sets of values of a polynomial. Throughout this paper, $X$ is an indeterminate, and
	\[
	f\in\Z[X]\quad\text{is a polynomial of degree }k\ge2\text{ with positive leading coefficient}.
	\]
	For positive integers $Q$ and $N$ we define
	\begin{equation}\label{eq:LS}
		\Delta_f(Q,N)=\sup\ \frac{\displaystyle\sum_{\substack{1\le q\le Q\\ f(q)\ge1}}\ \sum_{\substack{a=1\\ \gcd(a,f(q))=1}}^{f(q)}\Bigl|\sum_{M<n\le M+N}z_n\,\e\Bigl(\frac{an}{f(q)}\Bigr)\Bigr|^2}{\displaystyle\sum_{M<n\le M+N}|z_n|^2}\,,
	\end{equation}
	where the supremum is taken over all integers $M$ and all complex numbers $z_{M+1},\dots,z_{M+N}$ which are not all zero. The supremum is finite: by the Cauchy--Schwarz inequality each of the squares in the numerator is at most $N\sum|z_n|^2$. In other words, $\Delta_f(Q,N)$ is the best constant in the large sieve inequality for the moduli $f(q)$ with $q\le Q$. The problem is to determine its order of magnitude.
	
	\smallskip
	\emph{What to expect.} The numerator of \eqref{eq:LS} involves about $Q^{k+1}$ fractions $a/f(q)$. Taking a single non-zero coefficient $z_n$ shows that $\Delta_f(Q,N)$ is at least the number of these fractions, and taking $z_n=\e(-nx)$ for one fraction $x$ shows that $\Delta_f(Q,N)\ge N$. Hence
	\[
	\Delta_f(Q,N)\gg_fQ^{k+1}+N
	\]
	(see Section~\ref{sec:lower} for the details). If the fractions were evenly spaced, this would be the truth. For the powers $f=X^k$, that is, for the moduli $q^k$, Zhao \cite{Zh} conjectured that it is the truth up to a factor $Q^\eps$:
	\begin{equation}\label{eq:Zhao}
		\Delta_{X^k}(Q,N)\ll_\eps Q^\eps\bigl(Q^{k+1}+N\bigr)\qquad\text{for every }\eps>0.
	\end{equation}
	Halupczok and Munsch \cite[(2.7)]{HM} formulated the same expectation for general polynomials, also in several variables. In one variable their conjecture reads $\Delta_f(Q,N)\ll(QN)^{o(1)}(Q^{k+1}+N)$.
	
	Two upper bounds are immediate. For each $q$ the fractions $a/f(q)$ are $1/f(q)$-spaced, so by \eqref{eq:LSgeneral} they contribute at most $(N+f(q))\sum|z_n|^2$; summing over $q$ gives the first bound in
	\begin{equation}\label{eq:trivial}
		\Delta_f(Q,N)\ll_fQ^{k+1}+NQ\qquad\text{and}\qquad\Delta_f(Q,N)\ll_fN+Q^{2k}.
	\end{equation}
	The second bound holds because a value $f(q)$ is taken at most $k$ times, and distinct reduced fractions with denominators $\ll_fQ^k$ have distances $\gg_fQ^{-2k}$. So $\Delta_f(Q,N)$ has the order $Q^{k+1}$ for $N\le Q^k$ and the order $N$ for $N\ge Q^{2k}$, and the problem lies in the range $Q^k<N<Q^{2k}$. We call $N=Q^{k+1}$, where the two terms of \eqref{eq:Zhao} are equal, the \emph{critical point}. There the first bound in \eqref{eq:trivial} exceeds the expected bound by the factor $Q$.
	
	\smallskip
	\emph{Square moduli.} Zhao \cite{Zh} was the first to prove a non-trivial large sieve inequality for square moduli. His bound \cite[Thm.~2]{Zh} is
	\[
	\Delta_{X^2}(Q,N)\ll_\eps\log(2Q)\Bigl(Q^3+N^\eps\bigl(N\sqrt Q+Q^2\sqrt N\bigr)\Bigr).
	\]
	Improvements were given by Baier and Zhao in \cite{BZ05} and \cite{BZ}. The bound of \cite{BZ} is
	\begin{equation}\label{eq:BZ}
		\Delta_{X^2}(Q,N)\ll_\eps Q^\eps\Bigl(Q^3+N+\min\bigl(N\sqrt Q,\,Q^2\sqrt N\bigr)\Bigr).
	\end{equation}
	Its proof combines Weyl differencing, Poisson summation and an elementary counting argument, and it builds on the large sieve for sparse sets of moduli in \cite{B06}. At the critical point $N=Q^3$ the bound \eqref{eq:BZ} is larger than the conjectured one by the factor $Q^{1/2}$, and no unconditional bound with a smaller exponent than $\frac12$ at $N=Q^3$ has been known. Improvements were known only under hypotheses: under Hooley's Hypothesis $R^*$ for short Sali\'e sums the bound $Q^{1/3+\eps}N$ holds at $N=Q^3$ \cite{BR}, and savings also follow from hypotheses on additive energies of modular square roots \cite{BE}. In the other direction, Baier, Lynch and Zhao \cite{BLZ} showed that the factor $Q^\eps$ in \eqref{eq:Zhao} cannot be replaced by a constant:
	\[
	\Delta_{X^2}(Q,N)\gg Q^{\log2/((1+\eps)\log\log Q)}\bigl(Q^3+N\bigr)
	\]
	for infinitely many $Q$ and suitable $N\asymp Q^3$. An analogue of this lower bound for quadratic polynomials is contained in \cite[Thm.~1.3]{Co}.
	Baker \cite[Lemma~2]{Ba} proved a large sieve inequality for the moduli $gq^2$, by an argument which was suggested to him by Maynard; we come back to it in Section~\ref{sec:quad}.
	
	\smallskip
	\emph{Power moduli.} Let $k\ge3$, and put $\kappa=2^{k-1}$ and $\sigma=1/(2k(k-1))$. For the moduli $q^k$ the following bounds are known. Zhao \cite[Thm.~3]{Zh} and Baier and Zhao \cite[Thm.~1]{BZ05} proved
	\begin{align*}
		\Delta_{X^k}(Q,N)&\ll_\eps\log(2Q)\Bigl(Q^{k+1}+N^\eps\bigl(NQ^{1-1/\kappa}+N^{1-1/\kappa}Q^{1+k/\kappa}\bigr)\Bigr),\\
		\Delta_{X^k}(Q,N)&\ll_\eps(\log\log QN)^{k+1}\bigl(Q^{k+1}+N+N^{1/2+\eps}Q^k\bigr),
	\end{align*}
	respectively. Halupczok \cite[Thm.~1]{Ha12} proved
	\[
	\Delta_{X^k}(Q,N)\ll_\eps(QN)^\eps\bigl(Q^{k+1}+Q^{1-\sigma}N+Q^{1+k\sigma}N^{1-\sigma}\bigr).
	\]
	In the range $Q^k\le N\le Q^{2k}$, Munsch \cite[Thm.~1.2]{Mu} and Baier and Lynch \cite[Cor.~2.6]{BL} obtained
	\begin{equation}\label{eq:Mu}
		\Delta_{X^k}(Q,N)\ll_\eps(QN)^\eps Q^{1+\frac1{k+1}}N^{1-\frac1{k(k+1)}}\qquad\text{and}\qquad\Delta_{X^k}(Q,N)\ll_\eps Q^{\eps}Q^{1+\frac2{k+1}}N^{1-\frac2{k(k+1)}},
	\end{equation}
	respectively.
	
	\smallskip
	\emph{Polynomial moduli.} Halupczok \cite{Ha15} proved a large sieve inequality for polynomial moduli in several variables, by counting the fractions near a given one with multidimensional Weyl sums; the proof uses the mean value estimates of Parsell, Prendiville and Wooley \cite{PPW}. For a polynomial $f$ in one variable which has the coefficient $1$ at $X^j$ for some $2\le j\le k$, the bound \cite[Cor.~2.3]{Ha15} reads
	\begin{equation}\label{eq:Ha15}
		\Delta_f(Q,N)\ll_\eps Q^\eps\bigl(Q^{k+1}\log N+Q^{1-\sigma}N+Q^{1+(2k-j)\sigma}N^{1-\sigma}\bigr).
	\end{equation}
	For monic $f$ with $f(0)=0$ and $k\ge3$, Halupczok proved in \cite[Thm.~20]{Ha20}, in the range $Q^k\le N\le Q^{2k}$ and for the moduli $f(q)\gg Q^k$ with $Q\le q\le2Q$,
	\begin{equation}\label{eq:Ha20}
		\Delta_f(Q,N)\ll_\eps Q^\eps\Bigl(Q^{k+1}+\min\bigl(A_k,\,N^{1-\omega}Q^{1+(2k-1)\omega}\bigr)\Bigr),
	\end{equation}
	where
	\[
	A_k=NQ^{1-\frac1{k(k-1)}}+N^{1-\frac1{k(k-1)}}Q^{1+\frac1{k-1}}\qquad\text{and}\qquad\omega=\frac1{(k-1)(k-2)+2}.
	\]
	Halupczok and Munsch \cite[Thm.~2.3]{HM} proved the first bound in \eqref{eq:Mu} for polynomials with leading coefficient $\pm1$ and moduli of size $\gg Q^{k+o(1)}$; Munsch \cite[Rem.~1.4]{Mu} notes that his bound extends to polynomials of degree $k$. For every $f$ of degree $k\ge2$ and $Q^k\le N\le Q^{2k}$, Baker, Munsch and Shparlinski \cite[Thm.~1.5]{BMS} proved
	\begin{equation}\label{eq:BMS}
		\Delta_f(Q,N)\ll Q^{k+1}+Q^{o(1)}\bigl(NQ^{1-1/\kappa}+N^{1-1/\kappa}Q^{1+k/\kappa}\bigr),
	\end{equation}
	and McGrath \cite[Cor.~1.4]{McG} used the asymmetric additive energy of polynomials to prove, for $k\ge3$,
	\[
	\Delta_f(Q,N)\ll_\eps Q^\eps\bigl(NQ^{1/2}+N^{3/4}Q^{k/2+2-1/(50k)}\bigr).
	\]
	For quadratic polynomials the analogue of \eqref{eq:BZ} was stated by Corrigan in \cite{Co}, and in \cite{Co2} he used a generalised square sieve to establish the upper bound with the term $f(Q)\sqrt N$ alone in place of the minimum. A large sieve inequality for moduli of the form $4r-t^2$, where $r$ is a parameter, was proved by Shparlinski and Zhao \cite[Lemma~2.2]{SZ} in their work on isogeny classes of elliptic curves.
	
	\smallskip
	\emph{Function fields.} The analogous questions over polynomial rings over finite fields are the subject of a companion paper \cite{BBCL}. There Zhao's conjecture for power moduli can be proved in the limit in which the size of the finite field tends to infinity, and exact maxima of cluster sizes can be studied.
	
	\subsection{Main results}\label{sec:results}
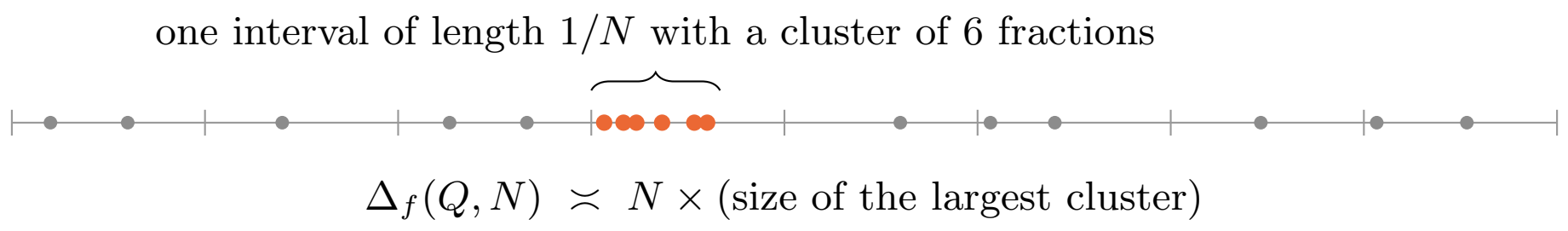
\begin{figure}[htb]
		\centering
		\begin{tikzpicture}[font=\footnotesize]
			\draw[black!40] (0,0) -- (12,0);
			\foreach \x in {0,1.5,...,12} \draw[black!40] (\x,-0.1) -- (\x,0.1);
			\foreach \x in {0.3,0.9,2.1,3.4,4.0,6.9,7.6,8.1,9.7,10.6,11.3} \fill[black!45] (\x,0) circle (1.5pt);
			\foreach \x in {4.6,4.75,4.85,5.05,5.3,5.4} \fill[hl] (\x,0) circle (1.8pt);
			\draw[decorate,decoration={brace,amplitude=4pt}] (4.5,0.25) -- node[above=4pt] {one interval of length $1/N$ with a cluster of $6$ fractions} (5.5,0.25);
			\node[below] at (6,-0.3) {$\Delta_f(Q,N)\ \asymp\ N\times(\text{size of the largest cluster})$};
		\end{tikzpicture}
		\caption{The principle behind all results of this paper.}\label{fig:cluster}
	\end{figure}

	All our results rest on one elementary principle (Lemma~\ref{lem:split}): \emph{up to an absolute constant, $\Delta_f(Q,N)$ is $N$ times the largest number of the fractions $a/f(q)$ in an interval of length $1/N$.} So upper and lower bounds for $\Delta_f(Q,N)$ are statements about \emph{clusters}, that is, about many fractions $a/f(q)$ in one short interval (Figure~\ref{fig:cluster}). On average an interval of length $1/N$ contains about $Q^{k+1}/N$ of the fractions; at the critical point this is a bounded number.
	
	For the precise statements we count the fractions in dyadic ranges of $q$. For real numbers $Y\ge1$ and $\delta>0$ let
	\begin{equation}\label{eq:Kdef}
		K_f(\delta,Y)=\sup_{\alpha\in\R}\ \#\Bigl\{(a,q)\in\Z^2:\ Y<q\le2Y,\ 1\le a\le f(q),\ \gcd(a,f(q))=1,\ \Bigl|\frac a{f(q)}-\alpha\Bigr|\le\delta\Bigr\}.
	\end{equation}
	This is the largest number of reduced fractions $a/f(q)$ with $Y<q\le2Y$ in an interval of length $2\delta$.
	
	\smallskip
	\emph{Clusters from ramification.} Our first result is a lower bound. It is governed by the following invariant. For a point $\xi\in\C$ let $e_f(\xi)$ be the multiplicity of $\xi$ as a zero of the polynomial $f(X)-f(\xi)$. This is the \emph{ramification index} of $f$ at $\xi$ \cite[\S3.1]{MZ}, and $e_f(\xi)\ge2$ exactly when $\xi$ is a critical point of $f$. We put
	\[
	\ell(f)=\max\bigl\{e_f(\xi):\ \xi\in\C,\ f(\xi)\ne0\bigr\}
	\]
	and call it the \emph{maximum ramification index of $f$ over the non-zero complex numbers}. We have $1\le\ell(f)\le k$. For example, $\ell(X^k)=1$, since the only critical point of $X^k$ is $0$ and the critical value is $0$; $\ell(X^2+1)=2$; and $\ell(X^k+c)=k$ for $c\ne0$. The polynomials with $\ell(f)=1$ form a small exceptional class: we shall see in Proposition~\ref{prop:class} that $\ell(f)=1$ if and only if $f$ belongs to
	\[
	\cE_k=\bigl\{c\,(X-\xi)^k:\ c,\xi\in\Q,\ c\ne0\bigr\}.
	\]
	
	\begin{maintheorem}[clusters from ramification]\label{main:lower}
		Let $\ell=\ell(f)$. There are constants $c>0$ and $C\ge1$, depending on $f$, such that for all $Q\ge C$ and $N\ge1$
		\begin{equation}\label{eq:mainlower}
			\Delta_f(Q,N)\ \ge\ c\,\Bigl(Q^{k+1}+N+Q^{1+k/\ell}N^{1-1/\ell}\Bigr).
		\end{equation}
		More precisely, if $\ell\ge2$, $Q\ge C$ and $CQ^k\le N\le cQ^{k+\ell}$, then some interval of length $1/(8N)$ contains at least
		\[
		c\,Q\,\Bigl(\frac{Q^k}N\Bigr)^{1/\ell}
		\]
		reduced fractions $a/f(q)$ with $Q/2<q\le Q$.
	\end{maintheorem}
	
	For $\ell=1$ the third term in \eqref{eq:mainlower} equals $Q^{k+1}$. For $\ell\ge2$ it is the largest of the three terms when $Q^k<N<Q^{k+\ell}$. At the critical point $N=Q^{k+1}$ it equals $Q^{1-1/\ell}N$, which gives the following.
	
	\begin{maincorollary}\label{main:Zhaofails}
		If $f\notin\cE_k$, then $\ell(f)\ge2$ and
		\[
		\Delta_f(Q,Q^{k+1})\gg_f Q^{1-1/\ell(f)}\,Q^{k+1}\ \ge\ Q^{1/2}\,Q^{k+1}.
		\]
	\end{maincorollary}
	
	So the analogue of Zhao's conjecture \eqref{eq:Zhao} for general polynomials, that is, the one-variable case of the conjecture of Halupczok and Munsch, holds at most for the polynomials in $\cE_k$. The simplest counterexample is $f=X^2+1$.
	
	\smallskip
	\emph{Upper bounds.} Our second result gives upper bounds. For quadratic polynomials $f=c_2X^2+c_1X+c_0$ we write $D=c_1^2-4c_2c_0$ for the discriminant.
	
	\begin{maintheorem}[upper bounds]\label{main:upper}
		Let $\eps>0$.
		\begin{enumerate}[label=\textup{(\roman*)}]
			\item For all $Q\ge1$ and $N\ge1$,
			\[
			\Delta_f(Q,N)\ll_{f,\eps}Q^\eps\bigl(Q^{k+1}+N+Q^k\sqrt N\bigr).
			\]
			\item Let $k=2$ and $D\ne0$. Then, for all sufficiently large $Q$ and all $N\ge1$,
			\[
			Q^3+N+Q^2\sqrt N\ \ll_f\ \Delta_f(Q,N)\ \ll_{f,\eps}\ Q^\eps\bigl(Q^3+N+Q^2\sqrt N\bigr).
			\]
			Moreover, for all sufficiently large $Y$ and $0<\delta\le\frac12$,
			\[
			1+Y^3\delta+Y^2\sqrt\delta\ \ll_f\ K_f(\delta,Y)\ \ll_{f,\eps}\ Y^\eps\bigl(1+Y^3\delta+Y^2\sqrt\delta\bigr).
			\]
			\item Let $k=2$ and $D=0$. Then, for all $Q\ge1$ and $N\ge1$,
			\[
			\Delta_f(Q,N)\ll_{f,\eps}Q^\eps\Bigl(Q^3+N+\min\bigl\{N\sqrt Q,\ Q^2\sqrt N,\ Q^{8/5}N^{3/5}\bigr\}\Bigr).
			\]
		\end{enumerate}
	\end{maintheorem}
	
	Part~(i) follows from a bound for clusters: we prove $K_f(\delta,Y)\ll_{f,\eps}Y^\eps(1+Y^{k+1}\delta+Y^k\sqrt\delta)$ for all sufficiently large $Y$ and $0<\delta\le1$ (Proposition~\ref{prop:uppergenr}). For quadratic polynomials the upper bound in (i) is due to Corrigan \cite[(1.3)]{Co2}. Our contribution in part~(ii) is the matching lower bound, which is the case $k=\ell=2$ of Theorem~\ref{main:lower}. Part~(ii) determines $\Delta_f(Q,N)$ up to the factor $Q^\eps$ for all quadratic polynomials of non-zero discriminant; in particular
	\[
	\Delta_{X^2+1}(Q,Q^3)=Q^{7/2+o(1)}.
	\]
	Part~(iii) says that for the quadratic polynomials in $\cE_2$, which are the polynomials $d(uX+v)^2$ with integers $d\ge1$, $u\ne0$ and $v$, the bound \eqref{eq:BZ} of Baier and Zhao and the bound of Theorem~\ref{main:squares} below both carry over from the squares.
	
	\smallskip
	\emph{Clusters from the leading monomial.} In degrees $k\ge3$ the clusters of Theorem~\ref{main:lower} are not always the largest ones, not even at the critical point.
	
	\begin{maintheorem}[clusters from the leading monomial]\label{main:leading}
		Let $f=X^k+g$ with a polynomial $g\in\Z[X]$ of degree $d<k$ such that $g(0)\ne0$. There are constants $c_f>0$ and $Q_f$ such that for all $Q\ge Q_f$ and every $N\ge1$
		\begin{equation}\label{eq:mainleading}
			\Delta_f(Q,N)\gg_f Q^{2-d/k}N^{1-1/k}.
		\end{equation}
		More precisely, if $Q\ge Q_f$ and $Q^{k-d}\le N\le c_fQ^{2k-d}$, then some interval of length $1/(8N)$ contains
		\[
		\gg_f\ Q\,\Bigl(\frac{Q^{k-d}}N\Bigr)^{1/k}
		\]
		reduced fractions $a/f(q)$ with $Q/2<q\le Q$.
	\end{maintheorem}
	
	The proof takes $q=bt$ and $a=t^k$. These fractions lie close to $1/b^k$, because $f(q)-q^k=g(q)$ is small compared with $f(q)$. Two examples show the effect. For $f=X^5+X+1$ we have $\ell(f)=2$. At the critical point $N=Q^6$, Theorem~\ref{main:leading} gives a cluster of size $\gg Q^{3/5}$, whereas Theorem~\ref{main:lower} gives $Q^{1/2}$. For $f=X^4+X+1$ and $N=Q^6$, Theorem~\ref{main:leading} still gives clusters of size $\gg Q^{1/4}$, whereas \eqref{eq:mainlower} gives only the order of a single fraction (Example~\ref{ex:leading}). So for $k\ge3$ the two constructions together do not determine the order of $\Delta_f(Q,N)$, which remains open. For the pure powers $X^k$, Zhao's conjecture \eqref{eq:Zhao} also remains open.
	
	\smallskip
	\emph{Square moduli.} For squares we obtain the following improvement on \eqref{eq:BZ}.
	
	\begin{maintheorem}[square moduli]\label{main:squares}
		For every $\eps>0$ and all $Q,N\ge1$,
		\[
		\Delta_{X^2}(Q,N)\ll_\eps Q^\eps\bigl(Q^3+N+Q^{8/5}N^{3/5}\bigr).
		\]
	\end{maintheorem}
	
\begin{figure}[htb]
		\centering
		\begin{tikzpicture}[x=5.2cm,y=6.2cm,font=\footnotesize]
			\draw[->] (1.95,0) -- (4.15,0) node[right] {$\nu$};
			\draw[->] (2,-0.02) -- (2,0.62) node[above] {excess exponent at $N=Q^\nu$};
			\foreach \x/\l in {2/2,2.5/{\frac52},3/3,3.5/{\frac72},4/4} \draw (\x,0.012) -- (\x,-0.012) node[below] {$\l$};
			\foreach \y/\l in {0.4/{\frac25},0.5/{\frac12}} \draw (2.015,\y) -- (1.985,\y) node[left] {$\l$};
			\draw[dashed,black!35] (3,0) -- (3,0.5);
			\draw[dashed,black!35] (2,0.5) -- (3,0.5);
			\draw[dashed,black!35] (2,0.4) -- (3,0.4);
			\draw[dashed,black!35] (2.75,0) -- (2.75,0.25);
			\node[below,black!60] at (2.75,-0.012) {$\frac{11}4$};
			\draw[thick,hl] (2,0) -- (3,0.5) -- (4,0);
			\draw[thick,bl] (2,0) -- (2.5,0) -- (3,0.5);
			\draw[very thick,gr] (2,0) -- (2.3333,0) -- (3,0.4) -- (4,0);
			\draw[thick,bl] (3,0.5) -- (4,0);
			\draw[thick,hl,dash pattern=on 4pt off 3pt] (3,0.5) -- (4,0);
			\node[hl,anchor=east] at (2.52,0.32) {$q^2+1$: the truth};
			\node[bl,anchor=west] at (3.35,0.4) {Baier--Zhao \eqref{eq:BZ}};
			\node[gr!70!black,anchor=west] at (3.07,0.2) {Theorem~\ref{main:squares}};
		\end{tikzpicture}
		\caption{Square moduli and the moduli $q^2+1$ at $N=Q^\nu$. Shown is the exponent $E$ for which the bound is $Q^{E+o(1)}(Q^3+N)$. Zhao's conjecture is $E=0$. For $q^2+1$ the true value of $E$ is known by Theorem~\ref{main:upper}(ii). For $\nu\ge3$ it coincides with the exponent in the bound of Baier and Zhao for $q^2$.}\label{fig:exponents}
	\end{figure}
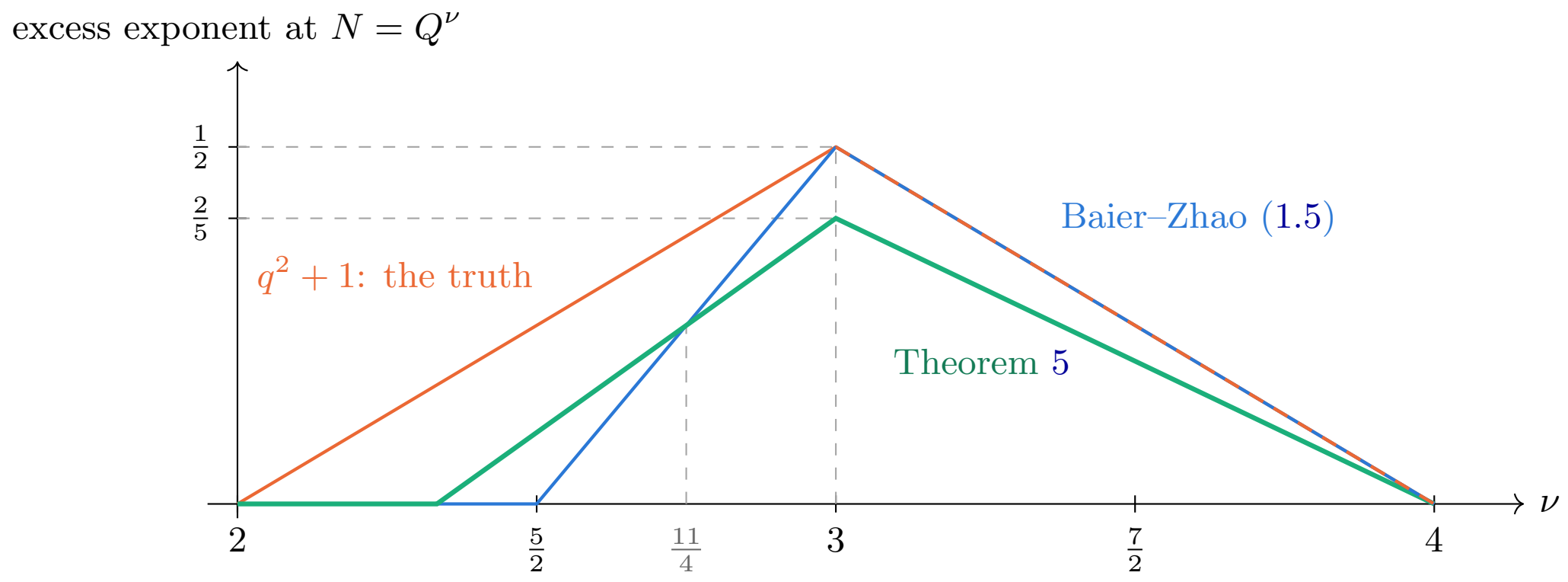

	At the critical point this reads $\Delta_{X^2}(Q,Q^3)\ll Q^{2/5+\eps}N$: an interval of length $Q^{-3}$ contains at most $Q^{2/5+\eps}$ reduced fractions $a/q^2$ with $q\le Q$. The theorem is better than \eqref{eq:BZ} exactly for $Q^{11/4}<N<Q^4$ (Figure~\ref{fig:exponents}). By Theorem~\ref{main:upper}(ii) it is false for the moduli $q^2+1$, so its proof has to use the arithmetic of squares. Zhao's conjecture corresponds to the exponent $0$ in place of $\frac25$, and Hypothesis $R^*$ gives $\frac13$ \cite{BR}.
	
	\subsection{Comparison with earlier results}\label{sec:lit}
	We now relate our results to the bounds listed in the first subsection.
	
	\smallskip
	\emph{The conjecture of Halupczok and Munsch.} Halupczok and Munsch \cite{HM} consider a polynomial $P\in\Z[X_1,\dots,X_m]$ of degree $k\ge2$ and the large sieve inequality for the moduli $P(\mathbf q)$ with $\mathbf q\in[Q,2Q)^m\cap\Z^m$, counted with multiplicity. On the basis of Zhao's heuristics they conjecture that the best constant is
	\[
	(QN)^{o(1)}\,r_P^*(Q)\,\bigl(Q^{m+k}+N\bigr),
	\]
	where $r_P^*(Q)$ is the largest number of points $\mathbf q$ of the box at which $P$ takes the same value \cite[(2.7)]{HM}. In one variable $r_P^*(Q)\le k$, so this is the bound $(QN)^{o(1)}(Q^{k+1}+N)$ in a dyadic range of $q$. Theorem~\ref{main:lower}, applied with $2Q-1$ in place of $Q$, gives $\gg_fQ^{1-1/\ell}$ reduced fractions $a/f(q)$ with $Q\le q<2Q$ in one interval of length $1/(8N)$ at $N=Q^{k+1}$. Hence the conjecture fails for every $f\notin\cE_k$. For the polynomials in $\cE_k$ our results do not decide it. For $f=X^k$ it is Zhao's conjecture in the range $Q^k\le N\le Q^{2k}$, and this is the essential case: a bound for the moduli $q^k$ carries over to all polynomials in $\cE_k$ (Remark~\ref{rem:Ek}).
	
	The one-variable conjecture of Halupczok \cite{Ha20}, which is quoted in \cite[\S2.2]{HM} in the form
	\[
	\Delta_f(Q,N)\ll(QN)^{o(1)}\bigl(Q^{k+1}+Q^{1+\frac1{k-1}}N^{1-\frac1{k(k-1)}}\bigr),
	\]
	is not affected, as we shall see in a moment. For $k=2$ its right-hand side is $Q^3+Q^2\sqrt N$, which by Theorem~\ref{main:upper}(ii) is the true order, up to $Q^{o(1)}$, for quadratic polynomials of non-zero discriminant and $Q^2\le N\le Q^4$.
	
	\smallskip
	\emph{The known upper bounds are compatible with Theorem~\ref{main:lower}.} Write $N=Q^\nu$ with $k\le\nu\le2k$, and recall that $\ell=\ell(f)\le k$. The cluster term $Q^{1+k/\ell}N^{1-1/\ell}$ of Theorem~\ref{main:lower} has the exponent $1+\nu-(\nu-k)/\ell$. We compare it with the exponents of the upper bounds.
	For the first bound in \eqref{eq:Mu}, which by \cite{HM} and \cite[Rem.~1.4]{Mu} holds for polynomial moduli, the exponent of the upper bound minus that of the cluster term is
	\[
	(\nu-k)\Bigl(\frac1\ell-\frac1{k(k+1)}\Bigr)\ \ge0.
	\]
	For the bound \eqref{eq:Ha20}, the corresponding differences for the second term of $A_k$ and for the term $N^{1-\omega}Q^{1+(2k-1)\omega}$ are
	\[
	(\nu-k)\Bigl(\frac1\ell-\frac1{k(k-1)}\Bigr)\ \ge0\qquad\text{and}\qquad\frac{\nu-k}\ell-(\nu-2k+1)\,\omega\ \ge\ (\nu-k)\Bigl(\frac1\ell-\omega\Bigr)\ \ge0,
	\]
	since $\omega\le1/k\le1/\ell$. So the minimum in \eqref{eq:Ha20} is at least the cluster term. The first of these two computations also applies to the conjecture of \cite{Ha20} displayed above. For the bound \eqref{eq:Ha15}, the exponent of $Q^{1+(2k-j)\sigma}N^{1-\sigma}$ minus that of the cluster term is
	\[
	(\nu-k)\Bigl(\frac1\ell-\sigma\Bigr)+(k-j)\,\sigma\ \ge0.
	\]
	For $k\ge3$ all these differences are positive when $\nu>k$. For $k=\ell=2$ the difference for the conjecture of \cite{Ha20} vanishes throughout $2\le\nu\le4$, in agreement with Theorem~\ref{main:upper}(ii). The bounds for the moduli $q^k$ are compatible trivially, since $\ell(X^k)=1$.
	
	\smallskip
	\emph{Quadratic moduli.} For quadratic polynomials $f$, the analogue of \eqref{eq:BZ} was stated in \cite[Thm.~1.1]{Co}, under size hypotheses which hold for a fixed $f$ and large $Q$:
	\begin{equation}\label{eq:Co}
		\Delta_f(Q,N)\ll_\eps f(Q)^\eps\Bigl(Qf(Q)+N+\min\bigl(N\sqrt Q,\,f(Q)\sqrt N\bigr)\Bigr).
	\end{equation}
	Theorem~\ref{main:upper}(ii) shows that this cannot hold when $D\ne0$: for a fixed $f$ and $0<\eta<\frac12$, the lower bound $Q^2\sqrt N$ exceeds the right-hand side of \eqref{eq:Co} by a power of $Q$ throughout the range $Q^{2+\eta}\le N\le Q^{3-\eta}$. In Section~\ref{sec:Co} we locate the step in \cite{Co} which has to be modified. The bound with $f(Q)\sqrt N$ in place of the minimum is correct; it was proved in \cite[(1.3)]{Co2}, and it is optimal by Theorem~\ref{main:upper}(ii). For $D=0$ the bound \eqref{eq:Co} holds for every $N$, with constants depending on $f$, by Theorem~\ref{main:upper}(iii).
	
	The moduli $q^2$ and $q^2+1$ behave differently below the critical point as well. For $f=X^2+1$, Theorem~\ref{main:upper}(ii) and \eqref{eq:trivial} give
	\[
	\Delta_f(Q,N)\ \begin{cases}\ \asymp Q^3&\text{for }N\le Q^2,\\ \ =Q^{2+o(1)}\sqrt N&\text{for }Q^2\le N\le Q^4,\\ \ \asymp N&\text{for }N\ge Q^4.\end{cases}
	\]
	On the other hand \eqref{eq:BZ} gives $\Delta_{X^2}(Q,N)\ll_\eps Q^\eps(Q^3+N\sqrt Q)$ for $Q^2\le N\le Q^3$. Since $Q^2\sqrt N/(Q^3+N\sqrt Q)\ge\frac12\min(N^{1/2}Q^{-1},\,Q^{3/2}N^{-1/2})$, we get
	\[
	\frac{\Delta_{X^2+1}(Q,N)}{\Delta_{X^2}(Q,N)}\gg_\eps Q^{-\eps}\min\bigl(N^{1/2}Q^{-1},\,Q^{3/2}N^{-1/2}\bigr)\qquad(Q^2\le N\le Q^3).
	\]
	For instance, at $N=Q^{5/2}$ the quotient is $\gg_\eps Q^{1/4-\eps}$.
	
	\smallskip
	\emph{Other upper bounds for polynomial moduli.} For $k\ge3$ and $Q^k<N<Q^{2k}$ the term $Q^k\sqrt N$ of Theorem~\ref{main:upper}(i) is smaller than the trivial term $NQ$ of \eqref{eq:trivial} precisely when $N>Q^{2k-2}$. So Theorem~\ref{main:upper}(i) is a result about the range near the spacing scale $N=Q^{2k}$: for $N=Q^{2k-\eta}$ with $0\le\eta\le2$, it determines $\Delta_f(Q,N)$ up to a factor $O_{f,\eps}(Q^{\eta/2+\eps})$. It complements the bounds obtained from Weyl sums. For example, for $k\ge3$ and $N=Q^\nu$ with $k<\nu<2k$, the three terms of Halupczok's bound \eqref{eq:Ha15} have the exponents $k+1$, $1-\sigma+\nu$ and $1+(2k-j)\sigma+\nu(1-\sigma)$, and the leading exponent in Theorem~\ref{main:upper}(i) is $k+\nu/2$. Solving the two non-constant inequalities shows that Halupczok's exponent is the smaller one for $\nu<\nu_j$ and the larger one for $\nu>\nu_j$, where
	\[
	\nu_j=2k-2+\min\Bigl\{2\sigma,\ \frac{2(j-2)\sigma}{1-2\sigma}\Bigr\}.
	\]
	For $k=j=2$, Theorem~\ref{main:upper}(i) has the smaller exponent throughout $2<\nu<4$.
	
	Finally we mention the approach through the distribution of the moduli in residue classes. For $S$ moduli in an interval of length $Q_0$, Baier \cite{B06} obtains a bound of the form $N+Q_0\,\Xi\,N^\eps(\sqrt N+S)$ under a hypothesis on the number of moduli in reduced residue classes in short intervals, in which $\Xi$ measures the deviation from uniform distribution. For polynomial moduli $S\asymp Q$ and $Q_0\asymp_fQ^k$. The clusters of Theorem~\ref{main:lower} come from residue classes in which the moduli are far from uniformly distributed: we shall find $\gg_fQ/b$ distinct moduli $f(q)$ in a single reduced residue class modulo $b^\ell$ (Lemma~\ref{lem:cong}). With $b^\ell\asymp Q$, a uniform estimate with the factor $\Xi$ would allow only $\ll\Xi(1+S/b^\ell)\ll_f\Xi$ moduli in this class; and for $N\asymp Q^{k+1}$, say, the modulus $b^\ell$ lies in the permitted range up to $\sqrt N$. Thus such a distribution estimate requires $\Xi\gg_fQ^{1-1/\ell}$ when $\ell\ge2$. The proof of Theorem~\ref{main:upper}(i) uses a different input, the number of representations of an integer as a difference $f(q_1)-f(q_2)$, and applies in the presence of these concentrations.
	
	\subsection{Outline of the methods}\label{sec:outline}
	We describe the three main ideas. In each case the task is to bound the number of fractions $a/f(q)$ in a short interval, from below or from above.
	
	\smallskip
	\emph{The clusters (Sections~\ref{sec:contact}--\ref{sec:conj}).} The simplest example is $f=X^2+1$. For integers $g\ge2$ and $A\ge1$ we have the identity
	\begin{equation}\label{eq:cluster}
		\frac{A^2}{(gA)^2+1}=\frac1{g^2}-\frac1{g^2\bigl((gA)^2+1\bigr)}.
	\end{equation}
	The fractions on the left-hand side are reduced, and those with $gA\asymp Q$ lie within $1/(g^2Q^2)$ of the single point $1/g^2$. With $g\asymp\sqrt Q$ we obtain about $\sqrt Q$ reduced fractions with denominators $q^2+1$, $q=gA\le Q$, in an interval of length about $Q^{-3}$. For $q^2$ in place of $q^2+1$ the same construction gives the fractions $A^2/(gA)^2$, which are not reduced and collapse to the single point $1/g^2$.
	
	What makes this work is that $f(q)\equiv1\pmod{g^2}$ for all $q\equiv0\pmod g$: the polynomial $f$ is constant modulo $g^2$ on a residue class modulo $g$, and the constant is a unit. This happens because $0$ is a critical point of $f$ with a non-zero critical value. In general, suppose that $f-f(\xi)$ has a zero of order $\ell\ge2$ at $\xi$ and that $f(\xi)\ne0$. If $\xi$ is an integer and $b$ is prime to $f(\xi)$, then $f(q)\equiv f(\xi)\pmod{b^\ell}$ for all $q\equiv\xi\pmod b$. Choose $y$ with $yf(\xi)\equiv1\pmod{b^\ell}$. Then for all these $q$ the number $a=(yf(q)-1)/b^\ell$ is an integer, and
	\[
	\frac a{f(q)}=\frac y{b^\ell}-\frac1{b^\ell f(q)}.
	\]
	All these fractions are reduced and lie within $1/(b^\ell Q^k)$ of the point $y/b^\ell$. There are about $Q/b$ of them, and choosing $b$ with $b^\ell Q^k\asymp N$ gives the cluster of Theorem~\ref{main:lower}. In general $\xi$ is an algebraic number, and the same works with $b$ a suitable value of the minimal polynomial of $\xi$. The condition $f(\xi)\ne0$ is what makes the fractions reduced; this is why the ramification over the value $0$ does not count in $\ell(f)$, and why the pure powers are exceptional.
	
	The clusters of Theorem~\ref{main:leading} are of a different kind. Here $f(q)$ is not constant modulo $b^k$ on the multiples $q=bt$ of $b$, but congruent to the number $g(q)$, which is small compared with $f(q)$.
	
	\smallskip
	\emph{Upper bounds for polynomial moduli (Section~\ref{sec:upperZ}).} Fix one reduced fraction $a_0/d_0$ in a cluster of fractions $a_i/d_i$ with denominators $d_i=f(q_i)\asymp Q^k$ in an interval of length $2\delta$. Each of the integers $a_0d_i-a_id_0$ is $O_f(\delta Q^{2k})$. Since $a_0$ is a unit modulo $d_0$, this restricts the differences $d_i-d_j$ to $O_f(1+\delta Q^{2k})$ possible values. Each non-zero value $n$ has few representations as $f(q_i)-f(q_j)$, since $q_i-q_j$ divides $n$. Counting the pairs $(i,j)$ gives Theorem~\ref{main:upper}(i).
	
	\smallskip
	\emph{The power saving for squares (Section~\ref{sec:below}).} Let $N=Q^3$, and consider the reduced fractions $a/q^2$ with $Q<q\le2Q$ in an interval of length $2/N$ around $\alpha$. The argument of Baier and Zhao approximates $\alpha$ by one fraction $v/u$. Then the integers $m=ua-vq^2$ are small, and
	\[
	q^2\equiv-\bar v\,m\pmod{|u|},
	\]
	where $\bar v$ is a multiplicative inverse of $v$ modulo $|u|$. So $q$ is a square root of a known residue, and for each $m$ there are $Q^{o(1)}$ fractions. But one has to pay for every value of $m$, and there are up to $\sqrt Q$ of them. This is the barrier $\frac12$; by Theorem~\ref{main:upper}(ii) it cannot be overcome for the moduli $q^2+1$.
	
	Our new idea is to admit only those $m$ which are squares, and in exchange to use many approximations $v/u$ instead of one. For parameters $T$ and $\delta=T^{-1/2+o(1)}$ we count in two ways the solutions of
	\begin{equation}\label{eq:introeq}
		ua-vq^2=t^2,\qquad T<t\le2T,\qquad|u|\le\delta q^2,
	\end{equation}
	where $a/q^2$ runs over the fractions in the interval. If $\gcd(q,u)=1$, the equation gives
	\[
	(t\bar q)^2\equiv-v\pmod{|u|},
	\]
	where $\bar q$ is a multiplicative inverse of $q$ modulo $|u|$. Now the square root is $t/q$, and its residue depends on the pair $(u,v)$ only, and no longer on $t$. So each pair $(u,v)$ has $Q^{o(1)}$ solutions, for all $t$ together. On the other hand, an average of quadratic Weyl sums shows that each fraction has at least $\frac12\delta T\approx\sqrt T$ solutions, apart from a small set of exceptional fractions which are counted separately. All pairs $(u,v)$ that occur lie in one parallelogram, which has the area $\ll Q^{o(1)}(Q/T+T^{3/2})$. Comparing the two counts, the number of fractions is
	\[
	\ll Q^{o(1)}\Bigl(\frac Q{T^{3/2}}+T\Bigr)=Q^{2/5+o(1)}\qquad\text{for }T=Q^{2/5}.
	\]
	A more detailed sketch, including the points which this description leaves out, is given in Section~\ref{sec:strategy}.
	
	\subsection{Notation and organisation}
	We write $A\ll B$ and $A=O(B)$ if $|A|\le CB$ with a constant $C$, which may depend on the quantities in the subscript. Further, $A\asymp B$ means $A\ll B\ll A$, and $Q^{o(1)}$ stands for a factor which is $\ll_\eps Q^\eps$ and $\gg_\eps Q^{-\eps}$ for every $\eps>0$. Constants depending on $f$ are allowed everywhere; we do not aim at uniformity in the coefficients of $f$. Table~\ref{tab:not} lists the notation that is used in the whole paper. We use the divisor bound $2^{\omega(s)}\le\tau(s)\ll_\eps s^\eps$ \cite[Thm.~315]{HW} without further comment.
	
	\begin{table}[htb]
		\centering\small
		\renewcommand{\arraystretch}{1.2}
		\begin{tabular}{p{3.2cm}p{11.4cm}}
			\toprule
			symbol & meaning\\
			\midrule
			$\e(x)$, $\|x\|$ & $\exp(2\pi ix)$; distance from $x$ to the nearest integer\\
			$f$, $k$, $c_k$ & the polynomial that generates the moduli; its degree; its leading coefficient, $c_k\ge1$\\
			$\Delta_f(Q,N)$ & the best constant in the large sieve inequality for the moduli $f(q)$, $q\le Q$; see \eqref{eq:LS}\\
			$K_f(\delta,Y)$ & the largest number of reduced fractions $a/f(q)$, $Y<q\le2Y$, in an interval of length $2\delta$; see \eqref{eq:Kdef}\\
			$f^{[j]}$ & the polynomial $f^{(j)}/j!$ (Section~\ref{sec:hasse})\\
			$e_f(\xi)$, $\ell(f)$ & the ramification index of $f$ at $\xi$; its maximum over the $\xi$ with $f(\xi)\ne0$\\
			$\cE_k$ & the polynomials $c(X-\xi)^k$ with $c,\xi\in\Q$, $c\ne0$; these are the $f$ with $\ell(f)=1$\\
			$\omega(s)$, $\tau(s)$, $\varphi(s)$ & number of distinct prime factors, number of divisors, Euler's function\\
			$\bar b$ & a multiplicative inverse of $b$ modulo the modulus in question\\
			$v_p(m)$, $p^e\,\|\,s$ & exponent of the prime $p$ in $m$; $p^e$ is the exact power of $p$ dividing $s$\\
			\bottomrule
		\end{tabular}
		\caption{General notation.}\label{tab:not}
	\end{table}
	
	Section~\ref{sec:prelim} collects the background: the reduction of the large sieve to the counting of fractions, Taylor expansions with integer coefficients, and square roots modulo an integer. Section~\ref{sec:contact} studies the invariant $\ell(f)$. Section~\ref{sec:lower} constructs the clusters from ramification and proves Theorem~\ref{main:lower}. Section~\ref{sec:conj} proves Theorem~\ref{main:leading}. Section~\ref{sec:upperZ} proves Theorem~\ref{main:upper} and discusses the bound stated in \cite{Co}. Section~\ref{sec:below} proves Theorem~\ref{main:squares}; it can be read independently of Sections~\ref{sec:contact}--\ref{sec:upperZ}.

	\section{Preliminaries}\label{sec:prelim}
	
	This section collects the background that is used throughout: the passage from the large sieve inequality to the counting of fractions in short intervals, Taylor expansions with integer coefficients, and two elementary lemmas on square roots modulo an integer.
	
	\subsection{From the large sieve to the counting of fractions}\label{sec:counting}
	The following lemma is the principle stated at the beginning of Section~\ref{sec:results}. It is well known; we include the short proof.
	
	\begin{lemma}[reduction to counting]\label{lem:split}
		Let $x_1,\dots,x_R$ be points of $\R/\Z$, not necessarily distinct, and let $N$ be a positive integer.
		\begin{enumerate}[label=\textup{(\roman*)}]
			\item Suppose that every closed arc of radius $\delta\in(0,\frac12]$ contains at most $K$ of the points $x_i$, counted with repetitions. Then for all integers $M$ and all complex numbers $z_n$
			\[
			\sum_{i=1}^R\Bigl|\sum_{M<n\le M+N}z_n\,\e(nx_i)\Bigr|^2\le K\bigl(N+\delta^{-1}\bigr)\sum_{M<n\le M+N}|z_n|^2 .
			\]
			\item Suppose that some arc of radius $1/(8N)$ contains $K$ of the points $x_i$. Then there are complex numbers $z_1,\dots,z_N$, not all zero, with
			\[
			\sum_{i=1}^R\Bigl|\sum_{n=1}^Nz_n\,\e(nx_i)\Bigr|^2\ge\frac{KN}2\sum_{n=1}^N|z_n|^2 .
			\]
		\end{enumerate}
	\end{lemma}
	
	\begin{proof}
		(i) Form a graph with the vertices $1,\dots,R$, in which $i\ne j$ are joined by an edge if $\|x_i-x_j\|<\delta$. All neighbours of $i$ belong to points in the arc of radius $\delta$ around $x_i$, so $i$ has at most $K-1$ neighbours. Hence a greedy colouring splits the vertices into $K$ classes without edges. The points of each class are $\delta$-spaced, and \eqref{eq:LSgeneral} applies to each class.
		
		(ii) Let $\alpha$ be the centre of the arc, and take $z_n=\e(-n\alpha)$. For each of the $K$ points $x_i$ in the arc we have $|2\pi n(x_i-\alpha)|\le\pi/4$ for $1\le n\le N$. So the real part of each term $z_n\e(nx_i)$ is at least $\cos(\pi/4)$, and $|\sum_nz_n\e(nx_i)|^2\ge N^2/2$, while $\sum_n|z_n|^2=N$.
	\end{proof}
	
	We apply the lemma to the fractions $a/f(q)$ in the definition \eqref{eq:LS} of $\Delta_f(Q,N)$. Part~(ii) gives the following: if some interval of length $1/(4N)$ contains $K$ of the fractions $a/f(q)$ with $q\le Q$, then
	\begin{equation}\label{eq:lowerprinciple}
		\Delta_f(Q,N)\ge\frac{KN}2 .
	\end{equation}
	For upper bounds we split the range of $q$ into dyadic blocks. Recall the definition \eqref{eq:Kdef} of $K_f(\delta,Y)$.
	
	\begin{lemma}[dyadic blocks]\label{lem:dyadic}
		Let $Y_0\ge1$ be such that $f$ is positive and strictly increasing on $[Y_0,\infty)$. There is a constant $C$, depending on $f$ and $Y_0$, such that for all integers $Q\ge1$ and $N\ge2$
		\[
		\Delta_f(Q,N)\le CN+4N\sum_{\substack{i\ge0\\ 2^iY_0<Q}}K_f\Bigl(\frac1N,\,2^iY_0\Bigr).
		\]
	\end{lemma}
	
	\begin{proof}
		For a single $q$ with $f(q)\ge1$ the fractions $a/f(q)$ are $1/f(q)$-spaced, so by \eqref{eq:LSgeneral} they contribute at most $(f(q)+N-1)\sum|z_n|^2\le f(q)N\sum|z_n|^2$ to the numerator in \eqref{eq:LS}. This accounts for the $q\le Y_0$, with $C=\sum f(q)$, the sum being over the $q\le Y_0$ with $f(q)\ge1$.
		
		Now let $Y=2^iY_0<Q$, and consider the fractions $a/f(q)$ with $Y<q\le\min(2Y,Q)$. They lie in $(0,1]$. A closed arc of radius $1/N\le\frac12$ of $\R/\Z$ meets $(0,1]$ in at most two intervals, each of length at most $2/N$, and each of these contains at most $K_f(1/N,Y)$ of the fractions. So Lemma~\ref{lem:split}(i) with $\delta=1/N$ and $K=2K_f(1/N,Y)$ shows that this block contributes at most $4N\,K_f(1/N,Y)\sum|z_n|^2$.
	\end{proof}
	
	The number of blocks is at most $1+\log_2Q$. So an upper bound for $K_f(1/N,Y)$ which is uniform in $Y\le Q$ gives, up to the factor $\log Q$, the bound $N$ times as large for $\Delta_f(Q,N)$. Together with \eqref{eq:lowerprinciple} this is the principle of Figure~\ref{fig:cluster}.
	
	\subsection{Taylor expansions with integer coefficients}\label{sec:hasse}
	For a polynomial $G$ with coefficients in a commutative ring $R$ we define the polynomials $G^{[j]}\in R[X]$ by
	\[
	G(X+Z)=\sum_{j\ge0}G^{[j]}(X)\,Z^j\qquad\text{in }R[X,Z].
	\]
	Thus $G^{[0]}=G$, $G^{[1]}=G'$, and $G^{[j]}=G^{(j)}/j!$ when $R$ contains $\Q$; these are the \emph{Hasse derivatives}. The point of the definition is that $G^{[j]}$ has coefficients in $R$: for $G\in\Z[X]$ the polynomial $G^{(j)}/j!$ has integer coefficients. For an element $\xi$ of a field containing $R$ we have
	\begin{equation}\label{eq:taylor}
		G(\xi+Z)-G(\xi)=\sum_{j\ge1}G^{[j]}(\xi)\,Z^j .
	\end{equation}
	So the multiplicity of $\xi$ as a zero of $G(X)-G(\xi)$ is the least $j\ge1$ with $G^{[j]}(\xi)\ne0$.
	
	\subsection{Square roots modulo an integer}\label{sec:sqrt}
	Square roots modulo an integer $s$ enter in Section~\ref{sec:upperZ} for quadratic moduli and in Section~\ref{sec:below} for square moduli. The following two lemmas are elementary.
	
	\begin{lemma}[square roots form few classes]\label{lem:sqrtZ}
		Let $s\ge1$ and $h$ be integers, and let $G=\gcd(h,s)$. The solutions of
		\[
		x^2\equiv h\pmod s
		\]
		form at most $2^{\omega(s)+1}$ residue classes modulo some divisor $s'$ of $s$ with $s'\ge s/\sqrt G$. Hence an interval $I$ contains at most $2^{\omega(s)+1}(|I|\sqrt G/s+1)$ such integers $x$.
	\end{lemma}
	
	\begin{proof}
		By the Chinese remainder theorem it suffices to show the following local statement. Let $p^e\,\|\,s$ and $p^j\,\|\,h$, with $j=\infty$ if $p^e\mid h$. Then the solutions modulo $p^e$ form at most $2$ classes ($4$ classes if $p=2$) modulo some $p^{e'}$ with $p^{e'}\ge p^{e-\min(j,e)/2}$.
		
		If $j\ge e$, the solutions are the $x\equiv0\pmod{p^{\lceil e/2\rceil}}$. This is one class, and $p^{\lceil e/2\rceil}\ge p^{e/2}$.
		
		If $j<e$, every solution has $v_p(x)=j/2$; there is none if $j$ is odd. So $x=p^{j/2}y$ with
		\[
		y^2\equiv h/p^j\pmod{p^{e-j}}\qquad\text{and}\qquad p\nmid y .
		\]
		This fixes $y$ in at most $2$ (resp.\ $4$) classes modulo $p^{e-j}$, hence $x$ in at most $2$ (resp.\ $4$) classes modulo $p^{e-j/2}$.
		
		The product of the local counts is at most $2^{\omega(s)+1}$, and each class modulo $s'$ meets $I$ in at most $|I|/s'+1$ integers.
	\end{proof}
	
	\begin{lemma}[square roots of a unit]\label{lem:roots}
		Let $s\ge1$. A unit modulo $s$ has at most $2^{\omega(s)+1}$ square roots modulo $s$.
	\end{lemma}
	
	\begin{proof}
		The square roots of a unit form a coset of the group of solutions of $x^2\equiv1\pmod s$, if there are any. Modulo an odd prime power this congruence has $2$ solutions, modulo a power of $2$ at most $4$.
	\end{proof}
	
	\section{Ramification over the non-zero values}\label{sec:contact}
	
	The clusters of Theorem~\ref{main:lower} come from points $\xi$ at which $f$ is ramified and does not vanish. In this section we study the invariant $\ell(f)$ that measures this ramification. We determine the polynomials for which it takes its extreme values $1$ and $k$, and we turn its definition into a divisibility statement between polynomials with rational coefficients (Lemma~\ref{lem:flatpoly}). This divisibility is all that the construction of the clusters in Section~\ref{sec:lower} uses. Let $\overline\Q\subset\C$ be the field of algebraic numbers.
	
	\begin{definition}\label{def:contact}
		Let $f\in\Q[X]$ have degree $k\ge1$. For $\xi\in\C$, the \emph{ramification index} $e_f(\xi)$ of $f$ at $\xi$ is the multiplicity of $\xi$ as a zero of the polynomial $f(X)-f(\xi)$ \cite[\S3.1]{MZ}. The \emph{maximum ramification index of $f$ over the non-zero complex numbers} is
		\[
		\ell(f)=\max\bigl\{e_f(\xi):\ \xi\in\C,\ f(\xi)\ne0\bigr\}.
		\]
	\end{definition}
	
	By \eqref{eq:taylor}, $e_f(\xi)$ is the least $j\ge1$ with $f^{(j)}(\xi)\ne0$. So $1\le e_f(\xi)\le k$, and $e_f(\xi)$ is $1$ plus the multiplicity of $\xi$ as a zero of $f'$. In particular $e_f(\xi)\ge2$ if and only if $\xi$ is a critical point of $f$; the number $f(\xi)$ is then the corresponding critical value. Critical points are zeros of $f'\in\Q[X]$, hence algebraic. So in the definition of $\ell(f)$ we may restrict $\xi$ to $\overline\Q$ as soon as $\ell(f)\ge2$. The maximum is over a non-empty set, because $f$ has only finitely many zeros.
	
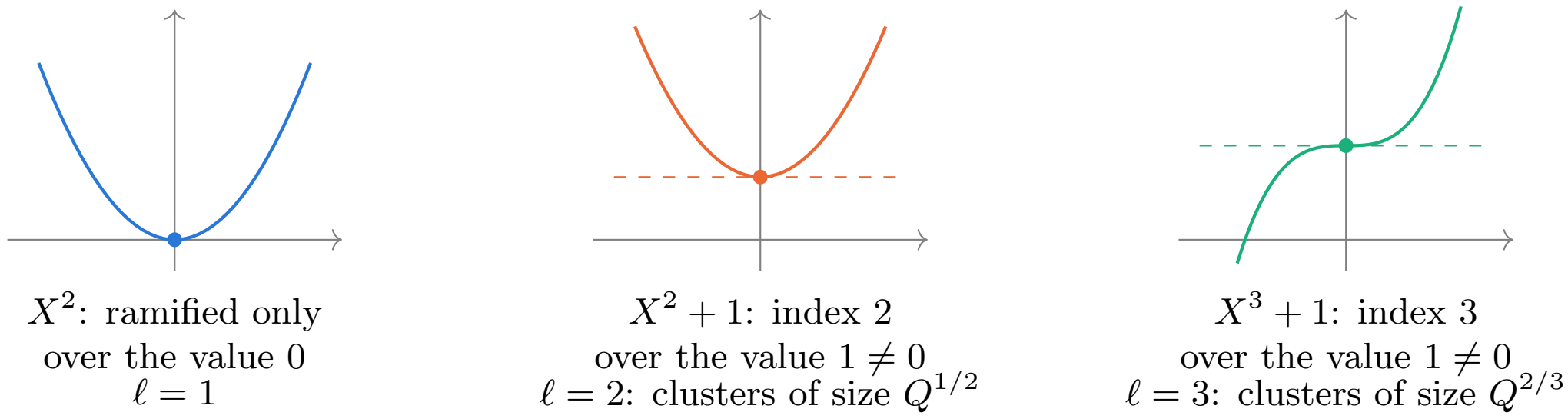
\begin{figure}[htb]
		\centering
		\begin{tikzpicture}[scale=0.9,font=\footnotesize]
			\begin{scope}
				\draw[->,black!50] (-1.6,0)--(1.6,0); \draw[->,black!50] (0,-0.3)--(0,2.2);
				\draw[thick,bl,domain=-1.3:1.3,samples=40] plot(\x,{\x*\x});
				\fill[bl] (0,0) circle (2pt);
				\node[below,align=center] at (0,-0.35) {$X^2$: ramified only\\ over the value $0$};
				\node at (0,-1.45) {$\ell=1$};
			\end{scope}
			\begin{scope}[xshift=5.6cm]
				\draw[->,black!50] (-1.6,0)--(1.6,0); \draw[->,black!50] (0,-0.3)--(0,2.2);
				\draw[thick,hl,domain=-1.2:1.2,samples=40] plot(\x,{\x*\x+0.6});
				\draw[dashed,hl] (-1.4,0.6)--(1.4,0.6);
				\fill[hl] (0,0.6) circle (2pt);
				\node[below,align=center] at (0,-0.35) {$X^2+1$: index $2$\\ over the value $1\neq0$};
				\node at (0,-1.45) {$\ell=2$: clusters of size $Q^{1/2}$};
			\end{scope}
			\begin{scope}[xshift=11.2cm]
				\draw[->,black!50] (-1.6,0)--(1.6,0); \draw[->,black!50] (0,-0.3)--(0,2.2);
				\draw[thick,gr,domain=-1.04:1.1,samples=40] plot(\x,{\x*\x*\x+0.9});
				\draw[dashed,gr] (-1.4,0.9)--(1.4,0.9);
				\fill[gr] (0,0.9) circle (2pt);
				\node[below,align=center] at (0,-0.35) {$X^3+1$: index $3$\\ over the value $1\neq0$};
				\node at (0,-1.45) {$\ell=3$: clusters of size $Q^{2/3}$};
			\end{scope}
		\end{tikzpicture}
		\caption{The invariant $\ell(f)$. The cluster sizes are those of Theorem~\ref{main:lower} at the critical point $N=Q^{k+1}$.}\label{fig:contactpic}
	\end{figure}

	The invariant $\ell(f)$ sees the ramification of $f$ in the fibres $f^{-1}(c)$ with $c\ne0$, but not in the fibre $f^{-1}(0)$ (Figure~\ref{fig:contactpic}).
	
	Recall that $\cE_k$ denotes the set of the polynomials $c(X-\xi_0)^k$ with $c,\xi_0\in\Q$ and $c\ne0$.
	
	\begin{proposition}[the extreme cases]\label{prop:class}
		Let $f\in\Q[X]$ have degree $k\ge2$.
		\begin{enumerate}[label=\textup{(\alph*)}]
			\item $\ell(f)=1$ if and only if $f\in\cE_k$.
			\item $\ell(f)=k$ if and only if $f=c(X-\xi_0)^k+d$ with $c,\xi_0,d\in\Q$ and $d\ne0$.
		\end{enumerate}
	\end{proposition}
	
	\begin{proof}
		(a) Write $f=c\prod_{i=1}^t(X-\xi_i)^{m_i}$ with distinct $\xi_i\in\C$. Then
		\[
		f'=c\prod_{i=1}^t(X-\xi_i)^{m_i-1}\cdot H\qquad\text{with}\qquad H=\sum_{i=1}^tm_i\prod_{j\ne i}(X-\xi_j).
		\]
		The polynomial $H$ has degree $t-1$ and leading coefficient $\sum_im_i=k$, and $H(\xi_i)=m_i\prod_{j\ne i}(\xi_i-\xi_j)\ne0$. If $t\ge2$, then $H$ has a zero $\eta$, and $f'(\eta)=0\ne f(\eta)$; so $\ell(f)\ge2$. If $t=1$, then $f=c(X-\xi_1)^k$, and $\xi_1\in\Q$ because the coefficient of $X^{k-1}$ is $-kc\xi_1$. Moreover $f'(\eta)=kc(\eta-\xi_1)^{k-1}\ne0$ whenever $f(\eta)\ne0$, so $\ell(f)=1$.
		
		(b) By \eqref{eq:taylor}, $e_f(\xi)=k$ means that $f(\xi+Z)-f(\xi)=cZ^k$, that is, $f=c(X-\xi)^k+f(\xi)$. As in (a), $\xi\in\Q$, and then $f(\xi)\in\Q$.
	\end{proof}
	
	\begin{example}\label{ex:contact}
		(i) $\ell(X^k)=1$, and $\ell(X^k+c)=k$ for $c\ne0$.
		
		(ii) Let $f=c_2X^2+c_1X+c_0$ with discriminant $D=c_1^2-4c_2c_0$. The only critical point is $\xi=-c_1/(2c_2)$, and the critical value is $f(\xi)=-D/(4c_2)$. So $\ell(f)=2$ if $D\ne0$, and $\ell(f)=1$ if $D=0$.
		
		(iii) $\ell(X^3+X+1)=2$: the critical points are $\pm\sqrt{-1/3}$, they are simple zeros of $f'$, and $f$ does not vanish there.
		
		(iv) $\ell(X^2(X-1))=2$: the critical points are $0$ and $\frac23$. The critical value at $0$ is $0$, so this point does not count, but $f(\frac23)=-\frac4{27}\ne0$.
	\end{example}
	
	In general the point $\xi$ at which the maximum $\ell(f)$ is attained is not rational, as in Example~\ref{ex:contact}(iii). For the construction of the clusters we need a statement about integers. The next lemma is the first step: it replaces the point $\xi$ by its minimal polynomial.
	
	\begin{lemma}[the flat polynomial]\label{lem:flatpoly}
		Let $f\in\Q[X]$, and let $\xi\in\overline\Q$ satisfy $f(\xi)\ne0$ and $e_f(\xi)=\ell\ge2$. Let $G\in\Q[X]$ be the minimal polynomial of $\xi$. Then
		\begin{equation}\label{eq:flatG}
			\gcd(G,f)=1\qquad\text{and}\qquad G^{\ell-j}\ \big|\ f^{[j]}\quad\text{in }\Q[X]\qquad(1\le j<\ell).
		\end{equation}
	\end{lemma}
	
	We call a non-constant polynomial $G$ with the properties \eqref{eq:flatG} a \emph{flat polynomial of order $\ell$} for $f$. For the quadratic polynomial of Example~\ref{ex:contact}(ii) the flat polynomial is $G=2c_2X+c_1$, up to a constant factor; for $X^k+c$ it is $G=X$; for $X^3+X+1$ it is $G=3X^2+1$.
	
	\begin{proof}
		Since $G$ is irreducible and $f(\xi)\ne0$, $G$ does not divide $f$, so $\gcd(G,f)=1$.
		
		Let $1\le j<\ell$. The polynomial $f(X)-f(\xi)$ is divisible by $(X-\xi)^\ell$, so its $j$-th derivative $f^{(j)}$ is divisible by $(X-\xi)^{\ell-j}$. Let $\xi'$ be another root of $G$. There is a field automorphism $\sigma$ of $\overline\Q$ with $\sigma(\xi)=\xi'$. We apply it to the coefficients of the identity $f^{(j)}=(X-\xi)^{\ell-j}h$, where $h\in\overline\Q[X]$. Since $f^{(j)}\in\Q[X]$, we see that $(X-\xi')^{\ell-j}$ divides $f^{(j)}$ as well. The roots of the irreducible polynomial $G$ are distinct, so $G^{\ell-j}$ divides $f^{(j)}=j!\,f^{[j]}$ in $\overline\Q[X]$, hence in $\Q[X]$.
	\end{proof}

	\section{Clusters from ramification: proof of Theorem~\ref{main:lower}}\label{sec:lower}
	
	In this section we prove Theorem~\ref{main:lower} and Corollary~\ref{main:Zhaofails}. Recall the plan from Section~\ref{sec:outline}: we look for a residue class $s\bmod b$ on which $f$ is constant modulo $b^\ell$, with a value prime to $b$. The proof has three steps.
	\begin{enumerate}[label=\arabic*.,leftmargin=2em]
		\item Lemma~\ref{lem:cong} turns the divisibilities \eqref{eq:flatG} between polynomials into such a congruence between integers, for suitable integers $s$ and $b$, and deduces the identity \eqref{eq:clusterid} which places the fractions near one point.
		\item In Section~\ref{sec:clusters} we choose the size of $b$ in terms of $Q$ and $N$ and count the fractions. This gives the cluster.
		\item Lemma~\ref{lem:card} supplies the term $Q^{k+1}$ in \eqref{eq:mainlower}: there are $\gg Q^{k+1}$ fractions in all.
	\end{enumerate}
	Throughout the section, $c_k\ge1$ denotes the leading coefficient of $f$.
	
	\subsection{Flat residue classes}
	
	\begin{lemma}[flat residue classes]\label{lem:cong}
		Let $G\in\Z[X]$ be a flat polynomial of order $\ell\ge2$ for $f$, and let $s_0$ be an integer with $A:=G(s_0)\ne0$. There is a positive integer $\Delta_0$, depending on $f$ and $G$ only, such that the following holds for every integer $s$ with
		\[
		s\equiv s_0\pmod{A^2\Delta_0}\qquad\text{and}\qquad G(s)\ne0 .
		\]
		The number $b=G(s)/A$ is a non-zero integer, $\gcd(b,f(s))=1$, and
		\begin{equation}\label{eq:flat}
			f(s+bu)\equiv f(s)\pmod{b^\ell}\qquad\text{for all }u\in\Z .
		\end{equation}
		Consequently, if $y$ is an integer with $yf(s)\equiv1\pmod{b^\ell}$, then for every integer $q\equiv s\pmod b$ with $f(q)\ne0$ the number
		\[
		a_q:=\frac{yf(q)-1}{b^\ell}
		\]
		is an integer which is prime to $f(q)$, and
		\begin{equation}\label{eq:clusterid}
			\frac{a_q}{f(q)}=\frac y{b^\ell}-\frac1{b^\ell f(q)} .
		\end{equation}
	\end{lemma}
	
	So $f$ is constant modulo $b^\ell$ on the residue class $s\bmod b$; we call it a \emph{flat class}. By \eqref{eq:clusterid}, the fractions $a_q/f(q)$ with $q$ in this class all lie within $1/(|b|^\ell|f(q)|)$ of the single point $y/b^\ell$.
	
	\begin{proof}
		\emph{Two identities.} Since $\gcd(G,f)=1$ and $G^{\ell-j}\mid f^{[j]}$ in $\Q[X]$, clearing denominators gives a positive integer $\Delta_0$ and polynomials $U,V,H_j\in\Z[X]$ with
		\begin{equation}\label{eq:bezoutid}
			UG+Vf=\Delta_0\qquad\text{and}\qquad\Delta_0f^{[j]}=G^{\ell-j}H_j\quad(1\le j<\ell).
		\end{equation}
		
		\emph{The integer $b$.} Let $s\equiv s_0\pmod{A^2\Delta_0}$. Then $G(s)\equiv G(s_0)=A\pmod{A^2\Delta_0}$, so $G(s)=A+A^2\Delta_0t$ with an integer $t$. Hence $b=1+A\Delta_0t$ is an integer, and it is prime to $\Delta_0$. Evaluating the first identity of \eqref{eq:bezoutid} at $s$ gives
		\[
		Ab\,U(s)+V(s)f(s)=\Delta_0 .
		\]
		So every common prime factor of $b$ and $f(s)$ divides $\Delta_0$; hence $\gcd(b,f(s))=1$.
		
		\emph{The congruence.} The second identity of \eqref{eq:bezoutid} gives $b^{\ell-j}\mid\Delta_0f^{[j]}(s)$, and as $b$ is prime to $\Delta_0$, we get $b^{\ell-j}\mid f^{[j]}(s)$ for $1\le j<\ell$. Now, by \eqref{eq:taylor},
		\[
		f(s+bu)-f(s)=\sum_{j\ge1}f^{[j]}(s)\,b^ju^j ,
		\]
		and every term is divisible by $b^\ell$: for $j<\ell$ because $b^{\ell-j}\mid f^{[j]}(s)$, and for $j\ge\ell$ trivially. This is \eqref{eq:flat}.
		
		\emph{The fractions.} For $q\equiv s\pmod b$, the congruence \eqref{eq:flat} gives $yf(q)\equiv yf(s)\equiv1\pmod{b^\ell}$, so $a_q$ is an integer. The identity $b^\ell a_q-yf(q)=-1$ shows that $a_q$ and $f(q)$ are coprime, and dividing it by $b^\ell f(q)$ gives \eqref{eq:clusterid}.
	\end{proof}
	
	\begin{example}\label{ex:classes}
		(i) For $f=X^2+1$ we have $G=X$ and $f-X\cdot G=1$, so we may take $\Delta_0=1$, $s_0=1$ and $A=1$. Then $s$ is arbitrary, $b=s$, the flat class is $q\equiv0\pmod b$, and $f(q)\equiv1\pmod{b^2}$ on it. With $y=1$ and $b=g\ge2$ the fractions \eqref{eq:clusterid} are those of \eqref{eq:cluster}.
		
		(ii) For $f=X^k+c$ with $c\ne0$, again $G=X$: for every $b$ prime to $c$ we have $f(q)\equiv c\pmod{b^k}$ for $q\equiv0\pmod b$.
		
		(iii) For $f=c_2X^2+c_1X+c_0$ with discriminant $D\ne0$ and $G=2c_2X+c_1$ we have $G^2-4c_2f=D$, so $\Delta_0=|D|$ is admissible. On the class $q\equiv s\pmod b$ one has $2c_2q+c_1\equiv0\pmod b$, hence
		\[
		4c_2f(q)=(2c_2q+c_1)^2-D\equiv-D\pmod{b^2}.
		\]
		
		(iv) For $f=X^3+X+1$ and $G=3X^2+1$ one has
		\[
		(6X^2-9X+4)\,G-(18X-27)\,f=31\qquad\text{and}\qquad f'=G,
		\]
		so $\Delta_0=31$ is admissible, with $s_0=0$ and $A=1$. For $s\equiv0\pmod{31}$ and $b=3s^2+1$ the polynomial $f$ is constant modulo $b^2$ on the class $s\bmod b$, although its critical points $\pm\sqrt{-1/3}$ are not rational. For instance $s=31$ gives $b=2884$.
	\end{example}
	
	\subsection{The clusters}\label{sec:clusters}
	We now prove the second assertion of Theorem~\ref{main:lower}. Let $\ell=\ell(f)\ge2$. We have to find constants $c,C>0$ such that for $Q\ge C$ and $CQ^k\le N\le cQ^{k+\ell}$ some interval of length $1/(8N)$ contains at least $cQ(Q^k/N)^{1/\ell}$ reduced fractions $a/f(q)$ with $Q/2<q\le Q$.
	
	\smallskip
	\emph{The flat polynomial.} Let $\xi$ be a point with $f(\xi)\ne0$ and $e_f(\xi)=\ell$, and let $G$ be an integer multiple of its minimal polynomial that has integer coefficients and positive leading coefficient. By Lemma~\ref{lem:flatpoly}, $G$ is a flat polynomial of order $\ell$ for $f$. Choose $s_0$ so large that $G$ is positive and increasing on $[s_0,\infty)$, put $A=G(s_0)$, and let $\Delta_0$ be as in Lemma~\ref{lem:cong}.
	
	\smallskip
	\emph{Flat classes of every size.} For $t=0,1,2,\dots$ put
	\[
	s_t=s_0+A^2\Delta_0\,t\qquad\text{and}\qquad b_t=\frac{G(s_t)}A .
	\]
	By Lemma~\ref{lem:cong} the $b_t$ are positive integers, and they increase with $t$. Since $s_{t+1}/s_t\to1$ and $G$ is a polynomial, $b_{t+1}/b_t\to1$. So there is $t_0$ with $b_{t+1}\le2b_t$ for $t\ge t_0$. Consequently, for every real $B\ge\max(2,b_{t_0})$ there is a $t$ with
	\begin{equation}\label{eq:bsize}
		B\le b_t\le2B :
	\end{equation}
	take the least $t\ge t_0$ with $b_t\ge B$.
	
	\smallskip
	\emph{The cluster.} Let $C$ be so large that for $Q\ge C$ the polynomial $f$ is increasing on $[Q/2,\infty)$ and $f(q)\ge\max(2,c_kQ^k/2^{k+1})$ for $q>Q/2$. Put
	\[
	\Lambda=\frac{2^{k+4}N}{c_kQ^k}\qquad\text{and}\qquad B=\Lambda^{1/\ell}.
	\]
	If $N\ge CQ^k$ with $C$ large enough, then $B\ge\max(2,b_{t_0})$, and \eqref{eq:bsize} gives $s=s_t$ and $b=b_t$ with
	\[
	\Lambda^{1/\ell}\le b\le2\Lambda^{1/\ell}.
	\]
	Since $\gcd(b,f(s))=1$, there is an integer $y$ with $1\le y<b^\ell$ and $yf(s)\equiv1\pmod{b^\ell}$. Let $q\equiv s\pmod b$ with $Q/2<q\le Q$, and let $a=a_q$ be as in Lemma~\ref{lem:cong}. Then $a/f(q)$ is a reduced fraction. The integer $yf(q)-1\ge f(q)-1\ge1$ is divisible by $b^\ell$, so $a\ge1$, and $a<yf(q)/b^\ell<f(q)$. By \eqref{eq:clusterid},
	\[
	0<\frac y{b^\ell}-\frac a{f(q)}=\frac1{b^\ell f(q)}\le\frac{2^{k+1}}{\Lambda c_kQ^k}=\frac1{8N}.
	\]
	So all these fractions lie in one interval of length $1/(8N)$. They are distinct, because $f$ is injective on $(Q/2,Q]$.
	
	\smallskip
	\emph{Their number.} If $N\le cQ^{k+\ell}$ with $c=c_k/(2^{k+4}8^\ell)$, then $b\le2\Lambda^{1/\ell}\le Q/4$. So the interval $(Q/2,Q]$ contains at least
	\[
	\frac Q{2b}-1\ \ge\ \frac Q{4b}\ \ge\ \frac Q{8\Lambda^{1/\ell}}\ \gg_f\ Q\,\Bigl(\frac{Q^k}N\Bigr)^{1/\ell}
	\]
	integers $q\equiv s\pmod b$. This proves the second assertion of Theorem~\ref{main:lower}. By \eqref{eq:lowerprinciple} it implies
	\begin{equation}\label{eq:clusterterm}
		\Delta_f(Q,N)\ \ge\ \frac N2\cdot\frac Q{8\Lambda^{1/\ell}}\ \gg_f\ Q^{1+k/\ell}N^{1-1/\ell}\qquad\text{for }CQ^k\le N\le cQ^{k+\ell}.
	\end{equation}
	
	\smallskip
	The size of $b$ is a free parameter: a larger $b$ gives fewer integers $q\le Q$ in the flat class, but a shorter interval. This is why clusters exist at every scale $Q^k\ll N\ll Q^{k+\ell}$, and why their size $Q(Q^k/N)^{1/\ell}$ decreases from $Q$ at $N\asymp Q^k$ to $1$ at $N\asymp Q^{k+\ell}$.
	
	\subsection{The number of fractions}
	For the term $Q^{k+1}$ in \eqref{eq:mainlower} we need to know that there are $\gg Q^{k+1}$ fractions. This is the content of the next lemma. A result of this kind is contained in an extension of \cite[Thm.~5.5A.2]{Sh}, which is proved by a different method; since we could not find this extension in the literature, we include a proof. For a polynomial $F\in\Z[X]$ and $d\ge1$ let $\varrho_F(d)$ be the number of $x$ modulo $d$ with $F(x)\equiv0\pmod d$.
	
	\begin{lemma}\label{lem:card}
		Let $f\in\Z[X]$ have degree $k\ge1$ and positive leading coefficient. There are $c>0$ and $Q_0$, depending on $f$, such that for $Q\ge Q_0$
		\[
		\sum_{Q/2<q\le Q}\varphi(f(q))\ \ge\ c\,Q^{k+1}.
		\]
	\end{lemma}
	
	\begin{proof}
		The idea is to write $\varphi(m)/m$ as a product over the prime factors of $m$, to treat the primes $p\le z$ by the sieve of Eratosthenes, and to show that the larger primes do not matter.
		
		Let $\kappa$ be the greatest common divisor of the coefficients of $f$. For a prime $p\nmid\kappa$ we have $\varrho_f(p)\le k$, and $\varrho_f(p)\le p$ always. So the product
		\[
		\mathfrak S=\prod_p\Bigl(1-\frac{\varrho_f(p)}{p^2}\Bigr)
		\]
		converges and is positive. Let $z\ge\max(2,\kappa)$ be a parameter and $Z$ the product of the primes $p\le z$. For an integer $m\ge2$,
		\begin{equation}\label{eq:phi}
			\frac{\varphi(m)}m=\prod_{p\mid m}\Bigl(1-\frac1p\Bigr)\ \ge\ \prod_{\substack{p\mid m\\ p\le z}}\Bigl(1-\frac1p\Bigr)-\sum_{\substack{p\mid m\\ p>z}}\frac1p,
		\end{equation}
		since $x\prod_i(1-x_i)\ge x-\sum_ix_i$ for $x,x_i\in[0,1]$. We sum this over $m=f(q)$ with $Q/2<q\le Q$, where $Q$ is so large that $2\le f(q)\le Q^{k+1}$ in this range.
		
		\emph{The small primes.} The first product in \eqref{eq:phi} equals $\sum_{d\mid\gcd(m,Z)}\mu(d)/d$. The interval $(Q/2,Q]$ contains $\varrho_f(d)\,(Q/(2d)+\theta_d)$ integers $q$ with $d\mid f(q)$, where $|\theta_d|\le1$. Hence
		\[
		\sum_{Q/2<q\le Q}\ \prod_{\substack{p\mid f(q)\\ p\le z}}\Bigl(1-\frac1p\Bigr)=\frac Q2\prod_{p\le z}\Bigl(1-\frac{\varrho_f(p)}{p^2}\Bigr)+O\bigl(2^{z}\bigr)\ \ge\ \frac Q2\,\mathfrak S-O(2^z),
		\]
		because $\varrho_f(d)\le d$ and $Z$ has at most $2^z$ divisors.
		
		\emph{The large primes.} For a prime $p$ with $z<p\le Q$ there are at most $k(Q/(2p)+1)$ integers $q\in(Q/2,Q]$ with $p\mid f(q)$. A number $f(q)\le Q^{k+1}$ has at most $k+1$ prime factors $p>Q$. Therefore
		\[
		\sum_{Q/2<q\le Q}\ \sum_{\substack{p\mid f(q)\\ p>z}}\frac1p\ \le\ \sum_{z<p\le Q}\frac kp\Bigl(\frac Q{2p}+1\Bigr)+\frac{(k+1)Q}{Q}\ \le\ \frac{kQ}z+k\sum_{p\le Q}\frac1p+k+1 ,
		\]
		using $\sum_{p>z}p^{-2}\le2/z$. The sum over $p\le Q$ is $O(\log\log Q)$.
		
		\emph{Conclusion.} Choosing $z$ with $k/z\le\mathfrak S/8$ and then $Q$ large, we get
		\[
		\sum_{Q/2<q\le Q}\frac{\varphi(f(q))}{f(q)}\ \ge\ \frac14\,\mathfrak S\,Q .
		\]
		Since $f(q)\gg_fQ^k$ for $q>Q/2$, the lemma follows.
	\end{proof}
	
	\subsection{Completion of the proofs}
	
	\begin{proof}[Proof of Theorem~\ref{main:lower}]
		The second assertion was proved in Section~\ref{sec:clusters}. For \eqref{eq:mainlower} we treat the three terms in turn. Taking a single non-zero coefficient $z_n$ in \eqref{eq:LS} shows that $\Delta_f(Q,N)$ is at least the number of fractions, which is at least $\sum_{Q/2<q\le Q}\varphi(f(q))\gg_fQ^{k+1}$ by Lemma~\ref{lem:card}. Taking $z_n=\e(-nx)$ for one of the fractions $x$ shows that $\Delta_f(Q,N)\ge N$.
		
		It remains to consider the third term $Q^{1+k/\ell}N^{1-1/\ell}$. If $\ell=1$, it equals $Q^{k+1}$. Let $\ell\ge2$. For $CQ^k\le N\le cQ^{k+\ell}$ the bound is \eqref{eq:clusterterm}. Outside this range the third term is dominated by one of the first two: for $N<CQ^k$ we have
		\[
		Q^{1+k/\ell}N^{1-1/\ell}\le C^{1-1/\ell}Q^{k+1},
		\]
		and for $N>cQ^{k+\ell}$ we have
		\[
		Q^{1+k/\ell}N^{1-1/\ell}=N\,\Bigl(\frac{Q^{k+\ell}}N\Bigr)^{1/\ell}\le c^{-1/\ell}N.\qedhere
		\]
	\end{proof}
	
	\begin{proof}[Proof of Corollary~\ref{main:Zhaofails}]
		We have $\ell=\ell(f)\ge2$ by Proposition~\ref{prop:class}(a). For $N=Q^{k+1}$ the third term of \eqref{eq:mainlower} is $Q^{1-1/\ell}N$.
	\end{proof}
	
	The corollary can be stated quantitatively: for $f\notin\cE_k$ the bound $\Delta_f(Q,N)\ll_\eps(QN)^\eps(Q^{k+1}+N)$ fails for every $\eps<(1-1/\ell(f))/(k+2)$. Indeed, at $N=Q^{k+1}$ this bound is at most $2C_\eps Q^{(k+2)\eps}N$, which is smaller than $cQ^{1-1/\ell}N$ for large $Q$ when $(k+2)\eps<1-1/\ell$.
	
	\section{Clusters from the leading monomial: proof of Theorem~\ref{main:leading}}\label{sec:conj}
	
	In this section $f=X^k+g$, where $g\in\Z[X]$ has degree $d<k$ and $g(0)\ne0$. The idea is simple. On the multiples $q=bt$ of an integer $b$, the leading term $q^k$ of $f(q)$ is exactly $b^kt^k$. So the fraction with numerator $t^k$ and denominator $f(q)$ differs from $1/b^k$ only through the remaining terms $g(q)$, and these are small if $d$ is small:
	\begin{equation}\label{eq:leadid}
		\frac1{b^k}-\frac{t^k}{f(q)}=\frac{f(q)-b^kt^k}{b^kf(q)}=\frac{g(q)}{b^k f(q)}\qquad(q=bt).
	\end{equation}
	Since $f(bt)\equiv g(0)\pmod t$, these fractions are reduced as soon as $t$ is prime to $g(0)$; this is where the condition $g(0)\ne0$ is used. The size of the clusters depends on the degree $d$ of the remainder.
	
	\begin{proof}[Proof of Theorem~\ref{main:leading}]
		\emph{The parameters.} Take $Q_f$ so large that for $Q\ge Q_f$ the polynomial $f$ is increasing on $(Q/2,Q]$ and $f(q)\ge q^k/2$ there. Choose $C_f\ge1$ so that
		\[
		\frac{|g(q)|}{f(q)}\le C_fQ^{d-k}\qquad\text{for }Q/2<q\le Q.
		\]
		Let $Q\ge Q_f$ and
		\[
		Q^{k-d}\le N\le c_fQ^{2k-d},
		\]
		where the constant $c_f>0$ will be chosen in a moment. We take
		\[
		b=\Bigl\lceil\bigl(16\,C_f\,N\,Q^{d-k}\bigr)^{1/k}\Bigr\rceil .
		\]
		Then $b\ge2$, because $NQ^{d-k}\ge1$, and $b\le Q/(4|g(0)|)$ if $c_f$ is sufficiently small and $Q_f$ sufficiently large.
		
		\emph{The fractions.} We consider the multiples $q=bt$ of $b$ in $(Q/2,Q]$ with
		\[
		t\equiv1\pmod{|g(0)|},
		\]
		and for each of them we take $a=t^k$. Since $f(bt)\equiv g(0)\pmod t$ and $\gcd(t,g(0))=1$, we have $\gcd(a,f(q))=1$. Also $1\le t^k<f(q)$, and the denominators $f(q)$ are distinct. By \eqref{eq:leadid} and the choice of $b$,
		\[
		\Bigl|\frac1{b^k}-\frac a{f(q)}\Bigr|=\frac{|g(q)|}{b^kf(q)}\le\frac{C_fQ^{d-k}}{b^k}\le\frac1{16N}.
		\]
		So all these fractions lie in one interval of length $1/(8N)$.
		
		\emph{Their number.} The admissible $t$ are the integers $t\equiv1\pmod{|g(0)|}$ in an interval of length $Q/(2b)$. Their number is at least
		\[
		\frac Q{2b\,|g(0)|}-1\ \ge\ \frac Q{4b\,|g(0)|},
		\]
		and $Q/b\gg_fQ(Q^{k-d}/N)^{1/k}$. This is the second assertion of the theorem.
		
		\emph{The bound for $\Delta_f(Q,N)$.} By \eqref{eq:lowerprinciple}, the cluster gives \eqref{eq:mainleading} for $Q^{k-d}\le N\le c_fQ^{2k-d}$:
		\[
		\Delta_f(Q,N)\gg_fN\cdot Q\,\Bigl(\frac{Q^{k-d}}N\Bigr)^{1/k}=Q^{2-d/k}N^{1-1/k}.
		\]
		If $N<Q^{k-d}$, the right-hand side of \eqref{eq:mainleading} is $\ll Q^{k+1-d}$, and if $N>c_fQ^{2k-d}$, it is $\ll_fN$. In both cases \eqref{eq:mainleading} follows from $\Delta_f(Q,N)\gg_fQ^{k+1}+N$, which is part of Theorem~\ref{main:lower}.
	\end{proof}
	
	\begin{example}[quintic and quartic moduli]\label{ex:leading}
		Put $f_k=X^k+X+1$ for $k=4,5$. Every zero $\xi$ of $f_k'$ is simple, has modulus $k^{-1/(k-1)}<1$, and satisfies $f_k(\xi)=1+(1-1/k)\xi\ne0$. Hence $\ell(f_k)=2$ in both cases.
		
		For $k=5$, take $Q=L^5$, $N=Q^6$ and $b=4L^2$ with a large positive integer $L$. For each multiple $q=bt$ of $b$ in $(Q/2,Q]$, the identity \eqref{eq:leadid} gives
		\[
		0<\frac1{b^5}-\frac{t^5}{f_5(q)}=\frac{q+1}{b^5 f_5(q)}\le\frac{2}{b^5q^4}\le\frac{32}{4^5 Q^6}=\frac1{32N}.
		\]
		These $Q/(2b)+O(1)\asymp Q^{3/5}$ reduced fractions give $\Delta_{f_5}(Q,Q^6)\gg NQ^{3/5}$ at the critical point. The three terms in \eqref{eq:mainlower} together have the order $NQ^{1/2}$ there.
		
		For $k=4$, take $Q=L^4$, $N=Q^6$ and $b=4L^3$. Now the same identity gives
		\[
		0<\frac1{b^4}-\frac{t^4}{f_4(q)}\le\frac{2}{b^4q^3}\le\frac{16}{4^4Q^6}=\frac1{16N}.
		\]
		There are $\asymp Q^{1/4}$ reduced fractions in this interval. Thus $\Delta_{f_4}(Q,Q^6)\gg NQ^{1/4}$ even at $N=Q^{k+\ell(f_4)}$, where the three terms in \eqref{eq:mainlower} have the total order $N$.
	\end{example}
	
	These examples show that the expression $Q^{k+1}+N+Q^{1+k/\ell(f)}N^{1-1/\ell(f)}$ of Theorem~\ref{main:lower} is not an upper bound for $\Delta_f(Q,N)$, up to factors $Q^{o(1)}$, for all $f$. For quadratic polynomials of non-zero discriminant it is, by Theorem~\ref{main:upper}(ii).

	\section{Upper bounds for polynomial moduli: proof of Theorem~\ref{main:upper}}\label{sec:upperZ}
	
	We now turn to upper bounds. By Lemma~\ref{lem:dyadic} it suffices to bound the number $K_f(\delta,Y)$ of reduced fractions $a/f(q)$ with $Y<q\le2Y$ in an interval of length $2\delta$; see \eqref{eq:Kdef}. The section is organised as follows.
	\begin{itemize}[leftmargin=2em]
		\item Section~\ref{sec:genpoly} proves the bound for $K_f(\delta,Y)$ for polynomials of arbitrary degree and deduces Theorem~\ref{main:upper}(i). The idea is that, in a cluster, fixing one fraction restricts all differences of denominators to few residue classes, and that a non-zero integer has few representations as a difference $f(q_i)-f(q_j)$.
		\item Section~\ref{sec:congf} gives a second proof for quadratic polynomials. It is not needed for Theorem~\ref{main:upper}, but it shows where the clusters sit: it records how the count depends on a rational approximation $b/s$ to the centre of the interval.
		\item Section~\ref{sec:quad} proves parts (ii) and (iii) of Theorem~\ref{main:upper}.
		\item Section~\ref{sec:Co} compares the result with the bound stated in \cite{Co}.
	\end{itemize}
	
	\subsection{Differences of polynomial values}\label{sec:genpoly}
	
	\begin{proposition}\label{prop:uppergenr}
		For all sufficiently large $Y$, depending on $f$, all $0<\delta\le1$ and every $\eps>0$,
		\[
		K_f(\delta,Y)\ll_{f,\eps}Y^\eps\bigl(1+Y^{k+1}\delta+Y^k\sqrt\delta\bigr).
		\]
	\end{proposition}
	
	\begin{proof}
		Let $Y$ be so large that on $(Y,2Y]$ the polynomial $f$ is strictly increasing and satisfies
		\[
		c\,Y^k\le f(q)\le C\,Y^k
		\]
		with positive constants $c,C$ depending only on $f$. Fix an interval $I$ of length $2\delta$. We distinguish two cases, according to whether one denominator can have several numerators in $I$ or not.
		
		\emph{Case 1: $2\delta CY^k\ge1$.} Each denominator $d=f(q)$ contributes at most $1+2\delta d\ll_f\delta Y^k$ numerators to $I$. Summing over $q$ gives $O_f(Y^{k+1}\delta)$ fractions.
		
		\emph{Case 2: $2\delta CY^k<1$.} Now each denominator occurs at most once in $I$. Let $a_i/d_i$ with $d_i=f(q_i)$ be the reduced fractions in $I$, let $K$ be their number, and suppose $K\ge2$. We fix one of them and denote it by $a_0/d_0$. The integers
		\[
		m_i=a_0d_i-a_id_0
		\]
		satisfy $|m_i|\le2\delta d_id_0\ll_f\delta Y^{2k}$. For every ordered pair $i\ne j$ the difference $n=d_i-d_j$ is non-zero, and
		\begin{equation}\label{eq:diffcong}
			|n|\ll_fY^k,\qquad a_0\,n\equiv m_i-m_j\pmod{d_0}.
		\end{equation}
		There are $O_f(1+\delta Y^{2k})$ possible integers $m_i-m_j$. Since $\gcd(a_0,d_0)=1$, each of them fixes the residue class of $n$ modulo $d_0$, and as $d_0\ge cY^k$, each class contains $O_f(1)$ integers $n$ with $|n|\ll_fY^k$. So the differences $n$ take $O_f(1+\delta Y^{2k})$ values.
		
		It remains to count the pairs $(i,j)$ with a given difference $n\ne0$. The non-zero integer $h=q_i-q_j$ divides $n=f(q_j+h)-f(q_j)$. For each of the at most $2\tau(|n|)$ signed divisors $h$ of $n$, the polynomial $f(X+h)-f(X)-n$ has degree $k-1$ and therefore at most $k-1$ integer roots $q_j$. By the divisor bound, each $n$ arises from $\ll_{f,\eps}Y^\eps$ pairs. Altogether
		\[
		K(K-1)\ll_{f,\eps}Y^\eps\bigl(1+\delta Y^{2k}\bigr),
		\]
		and taking square roots gives $K\ll_{f,\eps}Y^\eps(1+Y^k\sqrt\delta)$.
	\end{proof}
	
	\begin{proof}[Proof of Theorem~\ref{main:upper}\,(i)]
		For $N=1$ the bound follows from \eqref{eq:trivial}. Let $N\ge2$, and let $Y_0$ be so large that $f$ is positive and increasing on $[Y_0,\infty)$ and that Proposition~\ref{prop:uppergenr} holds for $Y\ge Y_0$. By Lemma~\ref{lem:dyadic} and Proposition~\ref{prop:uppergenr} with $\delta=1/N$,
		\[
		\Delta_f(Q,N)\ll_fN+N\sum_{\substack{i\ge0\\ Y=2^iY_0<Q}}K_f\Bigl(\frac1N,Y\Bigr)\ll_{f,\eps}N+N\log(2Q)\,Q^{\eps}\Bigl(1+\frac{Q^{k+1}}N+\frac{Q^k}{\sqrt N}\Bigr).
		\]
		This is $\ll_{f,\eps}Q^{2\eps}(Q^{k+1}+N+Q^k\sqrt N)$, and $\eps$ is arbitrary.
	\end{proof}
	
	\subsection{Quadratic polynomials: where the clusters sit}\label{sec:congf}
	In the rest of this section $k=2$. We write
	\[
	f=c_2X^2+c_1X+c_0\in\Z[X]\quad\text{with }c_2\ge1,\qquad D=c_1^2-4c_2c_0,\qquad Y_f=\frac{4(|c_1|+|c_0|)}{c_2}+1 .
	\]
	For $Y\ge Y_f$ and $Y<q_1<q_2\le2Y$ we have
	\begin{equation}\label{eq:fsize}
		\tfrac12c_2Y^2\le f(q_1)\le8c_2Y^2\qquad\text{and}\qquad f(q_2)-f(q_1)=(q_2-q_1)\bigl(c_2(q_1+q_2)+c_1\bigr)\ge c_2Y(q_2-q_1).
	\end{equation}
	So $f$ is injective on $(Y,2Y]$, and the reduced fractions $a/f(q)$ with $Y<q\le2Y$ are distinct.
	
	For quadratic polynomials a congruence $f(q)\equiv t\pmod s$ is a congruence for a square, and Lemma~\ref{lem:sqrtZ} counts its solutions.
	
	\begin{lemma}\label{lem:congf}
		For $s\ge1$, $t\in\Z$ and an interval $I$, the number of $q\in I$ with $f(q)\equiv t\pmod s$ is at most
		\[
		2^{\omega(4c_2s)+1}\Bigl(\frac{\sqrt{c_2}\,|I|\sqrt{\gcd(4c_2t+D,\,s)}}s+1\Bigr).
		\]
	\end{lemma}
	
	\begin{proof}
		Put $x=2c_2q+c_1$, so that $x^2=4c_2f(q)+D$. If $f(q)\equiv t\pmod s$, then
		\[
		x^2\equiv4c_2t+D\pmod{4c_2s}.
		\]
		The map $q\mapsto x$ is injective and sends $I$ into an interval of length $2c_2|I|$. We apply Lemma~\ref{lem:sqrtZ} with the modulus $4c_2s$, and we use $\gcd(4c_2t+D,4c_2s)\le4c_2\gcd(4c_2t+D,s)$.
	\end{proof}
	
	The next proposition bounds the number of fractions near a point $\alpha$ in terms of the denominator $s$ of a rational approximation to $\alpha$.
	
	\begin{proposition}\label{prop:upperZ}
		For $Y\ge Y_f$, $0<\delta\le1$ and $\alpha\in\R$ there is an integer $s$ with $1\le s\le\delta^{-1/2}$ such that the number of reduced fractions $a/f(q)$ with $Y<q\le2Y$ and $|a/f(q)-\alpha|\le\delta$ is at most
		\[
		1+C_f\,2^{\omega(4c_2s)}\tau(s)\bigl(Y^3\delta+Y^2\sqrt\delta\bigr),
		\]
		where $C_f$ depends only on $c_2$. If $\delta<1/(128c_2^2Y^4)$, this number is at most $1$. Consequently
		\[
		K_f(\delta,Y)\ll_{f,\eps}Y^\eps\bigl(1+Y^3\delta+Y^2\sqrt\delta\bigr).
		\]
	\end{proposition}
	
	\begin{proof}
		\emph{Very short intervals.} The fractions in question are distinct reduced fractions with denominators at most $8c_2Y^2$, by \eqref{eq:fsize}. So they differ by at least $1/(64c_2^2Y^4)$, and if $\delta<1/(128c_2^2Y^4)$, at most one of them lies within $\delta$ of $\alpha$. This also gives the last assertion: otherwise $s\le\delta^{-1/2}\le12c_2Y^2$, so $2^{\omega(4c_2s)}\tau(s)\ll_{f,\eps}Y^\eps$.
		
		\emph{Encoding.} By Dirichlet's approximation theorem \cite[\S11.3]{HW} there are coprime integers $b,s$ with
		\[
		1\le s\le\delta^{-1/2}\qquad\text{and}\qquad\Bigl|\alpha-\frac bs\Bigr|\le\frac{\sqrt\delta}s .
		\]
		Put $z=\alpha-b/s$, so that $s|z|\le\sqrt\delta$. For a fraction $a/f(q)$ within $\delta$ of $\alpha$ we put
		\[
		m=as-bf(q).
		\]
		Then
		\begin{equation}\label{eq:enc}
			\Bigl|\frac m{sf(q)}-z\Bigr|\le\delta\qquad\text{and}\qquad m\equiv-bf(q)\pmod s,
		\end{equation}
		and the fraction is determined by the pair $(q,m)$. If $m=0$, then $a/f(q)=b/s$, which happens for at most one fraction. With an inverse $\bar b$ of $b$ modulo $s$, the congruence in \eqref{eq:enc} reads
		\[
		f(q)\equiv-\bar b\,m\pmod s .
		\]
		So the number of fractions is at most $1+\#\mathcal E$, where $\mathcal E$ is the set of the pairs $(q,m)$ with $Y<q\le2Y$, $m\ne0$, the inequality in \eqref{eq:enc}, and this congruence. By Lemma~\ref{lem:congf}, the number of $q$ for a given $m$ is governed by
		\[
		G_m=\gcd\bigl(D-4c_2\bar b\,m,\ s\bigr).
		\]
		
		\emph{A divisor sum.} For $d\mid s$, the $m$ with $d\mid D-4c_2\bar bm$ lie in at most one residue class modulo $d/\gcd(d,4c_2)$, because $\bar b$ is a unit modulo $d$. Hence, for $L\ge1$,
		\[
		\#\{m:\ 0<|m|\le L,\ d\mid D-4c_2\bar bm\}\le\frac{8c_2L}d+1 .
		\]
		Using $\sqrt{G_m}\le\sum_{d\mid G_m}\sqrt d$ we get
		\[
		\sum_{0<|m|\le L}\sqrt{G_m}\ \le\ \sum_{d\mid s}\sqrt d\,\Bigl(\frac{8c_2L}d+1\Bigr)\ \le\ \tau(s)\bigl(8c_2L+\sqrt s\bigr),
		\]
		and therefore, for every $W>0$,
		\begin{equation}\label{eq:divsum}
			\sum_{0<|m|\le L}\Bigl(\frac{\sqrt{c_2}\,W\sqrt{G_m}}s+1\Bigr)\le\sqrt{c_2}\,\tau(s)\Bigl(\frac{8c_2WL}s+\frac W{\sqrt s}\Bigr)+2L .
		\end{equation}
		We shall apply this with $L$ a bound for $|m|$ and $W$ a bound for the length of the interval in which $q$ lies for a given $m$. There are two cases.
		
		\emph{Case 1: $|z|\le2\delta$.} Then
		\[
		|m|\le sf(q)\,(|z|+\delta)\le24c_2sY^2\delta=:L .
		\]
		If $L<1$, then $\mathcal E=\varnothing$. Otherwise, by Lemma~\ref{lem:congf} with $I=(Y,2Y]$ and by \eqref{eq:divsum} with $W=Y$,
		\[
		\#\mathcal E\le2^{\omega(4c_2s)+1}\Bigl(\sqrt{c_2}\,\tau(s)\Bigl(\frac{8c_2YL}s+\frac Y{\sqrt s}\Bigr)+2L\Bigr).
		\]
		We bound the three terms. First, $YL/s=24c_2Y^3\delta$. Second, $L\le24c_2Y^2\sqrt\delta$, because $s\delta\le\sqrt\delta$. Third, $L\ge1$ gives $s\ge1/(24c_2Y^2\delta)$, hence $Y/\sqrt s\le\sqrt{24c_2}\,Y^2\sqrt\delta$.
		
		\emph{Case 2: $|z|>2\delta$.} Then $m/(sf(q))$ lies in $[z-\delta,z+\delta]$, which does not contain $0$. So $m$ has the sign of $z$, and
		\[
		f(q)\in J_m:=\Bigl[\frac{|m|}{s(|z|+\delta)},\ \frac{|m|}{s(|z|-\delta)}\Bigr],\qquad|m|\le sf(q)(|z|+\delta)\le12c_2sY^2|z|=:L.
		\]
		As $z^2-\delta^2\ge\frac34z^2$, the interval $J_m$ has length at most
		\[
		\frac83\cdot\frac{|m|\delta}{sz^2}\le\frac{32c_2Y^2\delta}{|z|}.
		\]
		By \eqref{eq:fsize}, the $q$ with $f(q)\in J_m$ lie in an interval of length at most
		\[
		\frac{|J_m|}{c_2Y}\le W:=\frac{32Y\delta}{|z|}.
		\]
		If $L<1$, then $\mathcal E=\varnothing$. Otherwise Lemma~\ref{lem:congf} and \eqref{eq:divsum} give
		\[
		\#\mathcal E\le2^{\omega(4c_2s)+1}\Bigl(\sqrt{c_2}\,\tau(s)\Bigl(\frac{8c_2WL}s+\frac W{\sqrt s}\Bigr)+2L\Bigr).
		\]
		Again we bound the three terms. First, $WL/s\le384c_2Y^3\delta$. Second, $L\le12c_2Y^2\sqrt\delta$, as $s|z|\le\sqrt\delta$. For the middle term we use both constraints: $|z|>2\delta$ gives $W\le16Y$, and $L\ge1$ gives $|z|\ge1/(12c_2sY^2)$, hence $W\le384c_2Y^3\delta s$. Therefore
		\[
		\frac W{\sqrt s}\le\min\Bigl(\frac{16Y}{\sqrt s},\ 384c_2Y^3\delta\sqrt s\Bigr)\le\sqrt{6144c_2}\,Y^2\sqrt\delta.\qedhere
		\]
	\end{proof}
	
	\begin{remark}\label{rem:dequals}
		The proof shows which approximations $b/s$ produce clusters. The divisor $d=s$ in the sum leading to \eqref{eq:divsum} contributes when $s\mid D-4c_2\bar bm$. For $f=X^2+1$ and $s=g^2$ odd, this is $m\equiv-b\pmod s$, and the solutions satisfy $q\equiv0\pmod g$. For $W\ge\sqrt s$ an interval of length $W$ can therefore contain $\asymp W/\sqrt s$ such $q$. The fractions in \eqref{eq:cluster} realise this contribution. When $D=0$, the same divisor requires $s\mid4c_2m$, so every non-zero $m$ has $|m|\ge s/(4c_2)$.
	\end{remark}
	
	\subsection{\texorpdfstring{Proof of Theorem~\ref{main:upper}\,(ii) and (iii)}{}}\label{sec:quad}
	The upper bound in part~(ii) is the case $k=2$ of part~(i). As we mentioned in the introduction, it was proved by Corrigan \cite[(1.3)]{Co2} for quadratic polynomials which are positive and strictly increasing. Since we are not concerned with uniformity in the coefficients of $f$, this assumption can be removed: a polynomial with positive leading coefficient can be translated to meet it, and its finitely many earlier positive values contribute $O_f(N)$.
	
	\begin{proof}[Proof of Theorem~\ref{main:upper}\,(ii)]
		Let $D\ne0$. Then $\ell(f)=2$ by Example~\ref{ex:contact}(ii), and Theorem~\ref{main:lower} with $k=\ell=2$ gives $\Delta_f(Q,N)\gg_fQ^3+N+Q^2\sqrt N$. The upper bound is part~(i).
		
		It remains to prove the two bounds for $K_f(\delta,Y)$, where $0<\delta\le\frac12$. The upper bound is Proposition~\ref{prop:upperZ}. For the lower bound we treat the three terms $1$, $Y^3\delta$ and $Y^2\sqrt\delta$ separately.
		
		\emph{The term $1$.} We have $K_f(\delta,Y)\ge1$, because there is at least one fraction.
		
		\emph{The term $Y^3\delta$.} Let $F$ be the number of the reduced fractions $a/f(q)$ with $Y<q\le2Y$. Counting the pairs of a fraction $x$ and a point $\alpha\in[0,1)$ with $\|x-\alpha\|\le\delta$ gives
		\[
		\int_0^1\#\bigl\{x:\ \|x-\alpha\|\le\delta\bigr\}\,d\alpha=2\delta F.
		\]
		So some $\alpha$ has at least $2\delta F$ fractions within distance $\delta$ modulo $1$. They lie in at most three real intervals of length $2\delta$, so one of these contains at least $\frac23\delta F$ of them. Finally $F\gg_fY^3$ by Lemma~\ref{lem:card}.
		
		\emph{The term $Y^2\sqrt\delta$.} We apply the second assertion of Theorem~\ref{main:lower} with $Q=2Y$ and $N=\lceil1/(8\delta)\rceil$. It gives $\gg_fY^2\sqrt\delta$ fractions in an interval of length $\delta$ when $C_fY^{-4}\le\delta\le c_fY^{-2}$. Outside this range, $Y^2\sqrt\delta\ll_f1+Y^3\delta$.
	\end{proof}
	
	For part~(iii) we reduce the moduli $f(q)$ to square moduli. A reduction of this kind goes back to an argument of Maynard, which appears in the work of Baker \cite[Lemma~2]{Ba} on the moduli $gq^2$; our proof follows more or less the same idea.
	
	\begin{proof}[Proof of Theorem~\ref{main:upper}\,(iii)]
		\emph{The shape of $f$.} Let $D=0$. Then the double root of $f$ is rational. Choose a primitive linear polynomial $uX+v\in\Z[X]$ which vanishes there. Its square is primitive by Gauss's lemma, so comparing contents gives
		\[
		f=d\,(uX+v)^2\qquad\text{with an integer }d\ge1 .
		\]
		
		\emph{From the moduli $dj^2$ to the moduli $j^2$.} Let $j=|uq+v|$, and omit $j=0$, for which $f(q)=0$. Every reduced residue class modulo $dj^2$ lies in a reduced residue class $a_0\bmod j^2$. For each $a_0$ we enlarge the sum over the numerators to all $d$ lifts $a_0+hj^2$ with $0\le h<d$. Orthogonality of the additive characters modulo $d$ gives
		\[
		\sum_{h=0}^{d-1}\Bigl|\sum_nz_n\,\e\Bigl(\frac{n(a_0+hj^2)}{dj^2}\Bigr)\Bigr|^2
		=d\sum_{r\bmod d}\Bigl|\sum_{n\equiv r\,(d)}z_n\,\e\Bigl(\frac{na_0}{dj^2}\Bigr)\Bigr|^2 .
		\]
		In the inner sum on the right-hand side we write $n=dm+r$. This removes a constant phase of modulus one, and the remaining exponential is $\e(ma_0/j^2)$, where $m$ runs over an interval with at most $\lceil N/d\rceil$ integers. So for each $r$ we have a large sieve sum for the square moduli $j^2$. Each positive $j$ occurs for at most two values of $q$, and $j\le|u|Q+|v|$. Summing over $a_0$ and $j$, and applying the definition of $\Delta_{X^2}$ for each $r$, we obtain
		\begin{equation}\label{eq:quadratic-transfer}
			\Delta_f(Q,N)\le2d\,\Delta_{X^2}\bigl(|u|Q+|v|,\ \lceil N/d\rceil\bigr).
		\end{equation}
		
		\emph{Conclusion.} The bound \eqref{eq:BZ} of Baier and Zhao and Theorem~\ref{main:squares}, inserted on the right-hand side of \eqref{eq:quadratic-transfer}, give (iii).
	\end{proof}

	\begin{remark}[the exceptional class in every degree]\label{rem:Ek}
	The reduction \eqref{eq:quadratic-transfer} is not special to degree $2$. Let $k\ge2$, and let $f\in\cE_k$ have integer coefficients and positive leading coefficient. As in the first step of the proof,
	\[
	f=d\,(uX+v)^k
	\]
	with integers $d\ge1$, $u\ge1$ and $v$, where $uX+v$ is primitive. Put $j=|uq+v|$. Enlarging the sum over the numerators to the $d$ lifts $a_0+hj^k$ of a reduced residue class $a_0\bmod j^k$, and splitting $n=dm+r$ as before, we obtain in the same way
	\begin{equation}\label{eq:power-transfer}
	\Delta_f(Q,N)\le2d\,\Delta_{X^k}\bigl(uQ+|v|,\ \lceil N/d\rceil\bigr);
	\end{equation}
	here again each positive value of $j$ occurs for at most two values of $q$. So every upper bound for the pure powers $q^k$ carries over to all polynomials in $\cE_k$, with constants depending on $f$. In particular, Zhao's conjecture \eqref{eq:Zhao} for the moduli $q^k$ would settle the whole exceptional class of degree $k$.
	\end{remark}
	
	\subsection{The bound stated in \texorpdfstring{\cite{Co}}{[Co]}}\label{sec:Co}
	We explain the discrepancy between Theorem~\ref{main:upper}(ii) and the bound \eqref{eq:Co}, which was stated in \cite[Thm.~1.1]{Co}.
	
	\smallskip
	\emph{The two bounds are incompatible.} Take $f=X^2+1$ and $N=\lfloor Q^{5/2}\rfloor$. Theorem~\ref{main:upper}(ii) gives $\Delta_f(Q,N)\gg Q^{13/4}$. For the same polynomial the size hypotheses of \cite[Thm.~1.1]{Co}, namely $Q\ge|D|/c_2^2$ and $f(Q)\le3c_2Q^2$, hold for large $Q$, and \eqref{eq:Co} would give $\Delta_f(Q,N)\ll_\eps Q^{3+\eps}$. More generally, for a fixed strictly increasing quadratic polynomial of non-zero discriminant and $0<\eta<\frac12$, the lower bound of Theorem~\ref{main:upper}(ii) exceeds the bound in \eqref{eq:Co} by a power of $Q$ throughout $Q^{2+\eta}\le N\le Q^{3-\eta}$, if $\eps$ is sufficiently small in terms of $\eta$.
	
	\smallskip
	\emph{What remains true.} The bound with $f(Q)\sqrt N$ in place of the minimum is proved in \cite[(1.3)]{Co2}. Together with \eqref{eq:trivial}, it proves \eqref{eq:Co} with a constant depending on $f$ when $N\le f(Q)$ or $N\ge f(Q)^2/Q$. For $D=0$, the inequality \eqref{eq:quadratic-transfer} gives \eqref{eq:Co} for every $N$, with constants depending on $f$.
	
	\smallskip
	\emph{The step which has to be modified.} The last estimate in the proof of \cite[Lemma~4.3]{Co} bounds the number of integers in a short interval and in a prescribed residue class by the length of the interval divided by the modulus. But a single admissible integer can contribute even when this quotient tends to zero. We give the corresponding example on the grid used in \cite[Lemma~2.2]{Co}; we write $Y$ for the parameter which is called $M$ there.
	
	Take $f=X^2+1$, let $g$ be the least odd integer which is at least $\sqrt Y$, and put $r=g^2$ and $b=r-1$. In the notation of \cite{Co}, take $A=1$, $p'=p_r$ and $t=t_{rA}$, and set
	\[
	\delta=\frac1{4trY^2},\qquad
	h=\Bigl\lfloor\frac{4p'tY^2}{3r}+\frac12\Bigr\rfloor,\qquad
	z=\frac{p't}{hr^2}.
	\]
	Here $r\asymp Y$ and $p',t\ll_\eps Y^\eps$. We check that this centre satisfies the restrictions of \cite[Lemmas~2.2 and 3.5]{Co}. For large $Y$ we have $Y^{-4}\le\delta\le Y^{-2}$ and $r\le\delta^{-1/2}$. Moreover
	\[
	\frac h{p'/(r\sqrt\delta)}\asymp Y\sqrt{\frac tr}\to\infty,
	\]
	so the integer $h$ meets the lower cutoff of the grid; its coprimality condition with $A$ is automatic. Finally
	\[
	z=\frac{3/4+o(1)}{rY^2}\ \ge\ \frac{\delta t}{p'}.
	\]
	The bound asserted in \cite[Lemma~4.3]{Co} is $O_\eps(Y^\eps)$ here, since $Y^3\delta\ll1$ and $Y^2rz\ll1$.
	
	On the other hand, for every multiple $q$ of $g$ in $[1.04Y,1.35Y]$ we have $Y^2<f(q)<2Y^2$ and $f(q)\equiv1\pmod r$. The integer $a=(1+bf(q))/r$ satisfies $ar-bf(q)=1$. Consequently $a/f(q)$ is reduced, and for large $Y$
	\[
	\frac{a}{f(q)}-\frac br=\frac{1}{rf(q)}
	\in\frac{[0.54,0.93]}{rY^2}
	\subset[z-\delta t,\ z+\delta t].
	\]
	There are $\gg Y/g\gg\sqrt Y$ such fractions. They all have the same value $ar-bf(q)=1$.

	\smallskip
	\emph{A second point: the grid itself.} Independently of the example above, there is a difficulty in the passage from an arbitrary centre to the grid. It is asserted in \cite[(2.4)]{Co}, and used in \cite[Lemma~2.2]{Co}, that for every positive integer $A$ there is a number $t_A$, depending only on $A$, such that
	\[
	\sup_{R\ge1}\ \min_{\substack{r\le R\\ \gcd(r,2A)=1}}R\,\|r\alpha\|\le t_A
	\]
	for every irrational number $\alpha$. Such a bound cannot be uniform in $\alpha$. Take $A=1$, let $R>4t_1$, and let $\alpha$ be irrational with $|\alpha-\frac12|<1/(4R)$. Every admissible $r\le R$ is odd, so $r/2$ is half an odd integer, and $|r\alpha-r/2|<\frac14$. Hence $\|r\alpha\|>\frac14$ and
	\[
	R\,\|r\alpha\|>\frac R4>t_1 .
	\]
	The example above is not affected by this: it only uses the conditions on the size and the coprimality of its parameters. And the upper bound of Theorem~\ref{main:upper}(i) for quadratic polynomials does not use the grid at all; it follows directly from Proposition~\ref{prop:uppergenr}.

	\section{Square moduli: proof of Theorem~\ref{main:squares}}\label{sec:below}
	
	In this section $f=X^2$, so the fractions are the $a/q^2$, and we prove Theorem~\ref{main:squares}. The section can be read independently of Sections~\ref{sec:contact}--\ref{sec:upperZ}: it uses only the reduction to counting (Lemmas~\ref{lem:split} and~\ref{lem:dyadic}), the bound for the number of square roots of a unit (Lemma~\ref{lem:roots}), Dirichlet's approximation theorem, the divisor bound, and the trivial bound \eqref{eq:trivial} for $N=1$. For $1\le n\le Q^{O(1)}$ we use the divisor bound in the form $2^{\omega(n)}\le\tau(n)=Q^{o(1)}$.
	
	\subsection{The counting problem}
	For $Q\ge2$, $N>0$ and $\alpha\in\R$ let
	\[
	\mathcal S(\alpha)=\Bigl\{(a,q)\in\Z^2:\ Q<q\le2Q,\ \gcd(a,q)=1,\ \Bigl|\frac a{q^2}-\alpha\Bigr|\le\frac1N\Bigr\}.
	\]
	This is the set of the reduced fractions $a/q^2$ with $Q<q\le2Q$ in the interval of length $2/N$ around $\alpha$. Two quantities govern its size:
	\[
	\mu=\frac{Q^3}N\qquad\text{and}\qquad\kappa=\frac{Q^4}N .
	\]
	The first is the order of the average number of fractions in such an interval. The second measures the length of the interval against the spacing: distinct pairs $(a,q)$ give distinct reduced fractions $a/q^2$, and two of them differ by at least $(2Q)^{-4}$, so that
	\begin{equation}\label{eq:Strivial}
		\#\mathcal S(\alpha)\le1+32\kappa .
	\end{equation}
	At the critical point $N=Q^3$ we have $\mu=1$ and $\kappa=Q$. The following proposition is the heart of the matter. Theorem~\ref{main:squares} will follow from it by Lemma~\ref{lem:dyadic}.
	
	\begin{proposition}\label{thm:thin}
		For $\eps>0$, $Q\ge2$ and $Q^2\le N\le Q^4$,
		\[
		\sup_{\alpha\in\R}\#\mathcal S(\alpha)\ll_\eps Q^\eps\bigl(1+\mu+\kappa^{2/5}\bigr).
		\]
	\end{proposition}
	
	\subsection{The strategy}\label{sec:strategy}
	Before we go into the details, we describe the proof of Proposition~\ref{thm:thin}. One principle is used throughout:
	\begin{quote}
		\emph{a square root of a given residue modulo an integer, which lies in a short real interval, has few possibilities.}
	\end{quote}
	What changes, in comparison with earlier work, is the choice of the number that plays the role of the square root.
	
	\smallskip
	\emph{The count of Baier and Zhao.} Let $b/r$ be a rational approximation to $\alpha$ with $r\le\sqrt N$, as given by Dirichlet's theorem. For $(a,q)\in\mathcal S(\alpha)$ the integer $m=ra-bq^2$ is then $O(\sqrt\kappa)$, and
	\[
	b\,q^2\equiv-m\pmod r .
	\]
	So $q$ is a square root of a residue which is determined by $m$, and $q$ is confined to a short interval. Counting these square roots and summing over $m$, one finds that for every $B\ge1$ at most
	\[
	Q^{o(1)}\,(1+B+\mu)
	\]
	fractions have $0<|m|\le B$ (Lemma~\ref{lem:label}). Only the sum over $m$ is bounded in this way: a single value of $m$ can occur for many fractions when $q$ and $r$ have a large common factor, but such $m$ are multiples of a large integer and therefore rare. With $B\asymp\sqrt\kappa$ we obtain the bound $Q^{o(1)}(1+\mu+\sqrt\kappa)$, which corresponds to the term $Q^2\sqrt N$ in the bound \eqref{eq:BZ} of Baier and Zhao. Its weakness is that one has to pay for every value of $m$.
	
	\smallskip
	\emph{The new equation.} We admit only those values of $m$ which are squares $t^2$ with $T<t\le2T$, and in exchange we use all approximations $v/u$ in a certain range instead of a single one. We call a quintuple $(a,q,t,u,v)$ with $(a,q)\in\mathcal S(\alpha)$ a \emph{solution} if
	\[
	ua-vq^2=t^2,\qquad T<t\le2T,\qquad|u|\le\delta q^2 .
	\]
	Here $T$ and $\delta$ are parameters, which will be $T=\kappa^{2/5}$ and $\delta=T^{-1/2+\eta}$ with a small $\eta>0$. If $\gcd(q,u)=1$, the equation gives
	\[
	\Bigl(\frac tq\Bigr)^2\equiv-v\pmod{|u|} .
	\]
	Now the square root is $t/q$, and the residue depends on the pair $(u,v)$ only. The integer $t$, which replaces $m$, has disappeared from the right-hand side. We count the solutions in two ways.
	
	\smallskip
	\emph{Each fraction has many solutions (Section~\ref{sec:eq}).} Let $\bar a$ be the inverse of $a$ modulo $q^2$. For given $(a,q)$ and $t$, a solution exists if and only if $\|t^2\bar a/q^2\|\le\delta$, and then it is unique. So we need many $t\in(T,2T]$ for which $t^2\theta$, $\theta=\bar a/q^2$, is close to an integer. If the values $t^2\theta$ are equidistributed modulo $1$, there are about $\delta T$ such $t$. An average of quadratic Weyl sums (Section~\ref{sec:weyl}) shows that this is so, in the form of the lower bound $\frac12\delta T$, unless $\theta$ has a rational approximation with denominator at most $T$. The fractions with such a $\theta$ are called \emph{major}, the others \emph{minor}. On average over the frequency $h$, the Weyl sums are of size about $\sqrt T$ (Lemma~\ref{lem:weyl}), and this is why $\delta$ has to be a little larger than $T^{-1/2}$. The major fractions are few: they are counted by the argument of Baier and Zhao, applied to $O(1)$ pairs $(r,b)$ (Lemma~\ref{lem:major}).
	
	\smallskip
	\emph{Each pair $(u,v)$ has few solutions (Section~\ref{sec:multiplicity}).} By the principle above, applied to the square root $t/q$, a pair $(u,v)$ has at most
	\[
	Q^{o(1)}\Bigl(1+\frac{\mu\sqrt{\gcd(u,v)}}T\Bigr)
	\]
	solutions (Lemma~\ref{lem:mult}). The factor $\sqrt{\gcd(u,v)}$ appears because common factors of $u$ and $v$ have to be removed before the congruence can be used.
	
\begin{figure}[htb]
		\centering
		\begin{tikzpicture}[x=0.78cm,y=0.78cm,font=\footnotesize]
			\fill[bl!10] (7,3.39) -- (7,2.49) -- (-7,-3.39) -- (-7,-2.49) -- cycle;
			\draw[bl] (7,3.39) -- (7,2.49) -- (-7,-3.39) -- (-7,-2.49) -- cycle;
			\draw[->,black!45] (-7.8,0) -- (8,0) node[right] {$u$};
			\draw[->,black!45] (0,-3.7) -- (0,3.9) node[above] {$v$};
			\foreach \x in {-7,...,7} \foreach \y in {-3,...,3} \fill[black!25] (\x,\y) circle (0.9pt);
			\draw[hl,thick] (-5.6,-2.8) -- (5.6,2.8) node[above right=-1pt] {$\R\mathbf w_0$};
			\foreach \p in {(2,1),(4,2),(-2,-1),(-4,-2)} \fill[hl] \p circle (2.2pt);
			\foreach \p in {(1,0),(3,1),(5,2),(7,3),(-1,0),(-3,-1),(-5,-2),(-7,-3)} \fill[bl] \p circle (2.2pt);
			\fill (0,0) circle (2.2pt);
			\node[hl,above left] at (2,1) {$\mathbf w_0$};
			\node[bl,anchor=west] at (7.15,2.2) {$\Pi$};
			\draw[<->,black!60] (-7,-3.75) -- (7,-3.75) node[midway,below] {$|u|\le4Q^2\delta$};
		\end{tikzpicture}
		\caption{The parallelogram $\Pi$ of Lemma~\ref{lem:par}: a thin neighbourhood of the line $v=\alpha u$. The integer points off the line $\R\mathbf w_0$ through the shortest point $\mathbf w_0$ are $O(A)$ in number (Lemma~\ref{lem:offline}); the points $i\mathbf w_0$ on the line are treated separately.}\label{fig:Pi}
	\end{figure}
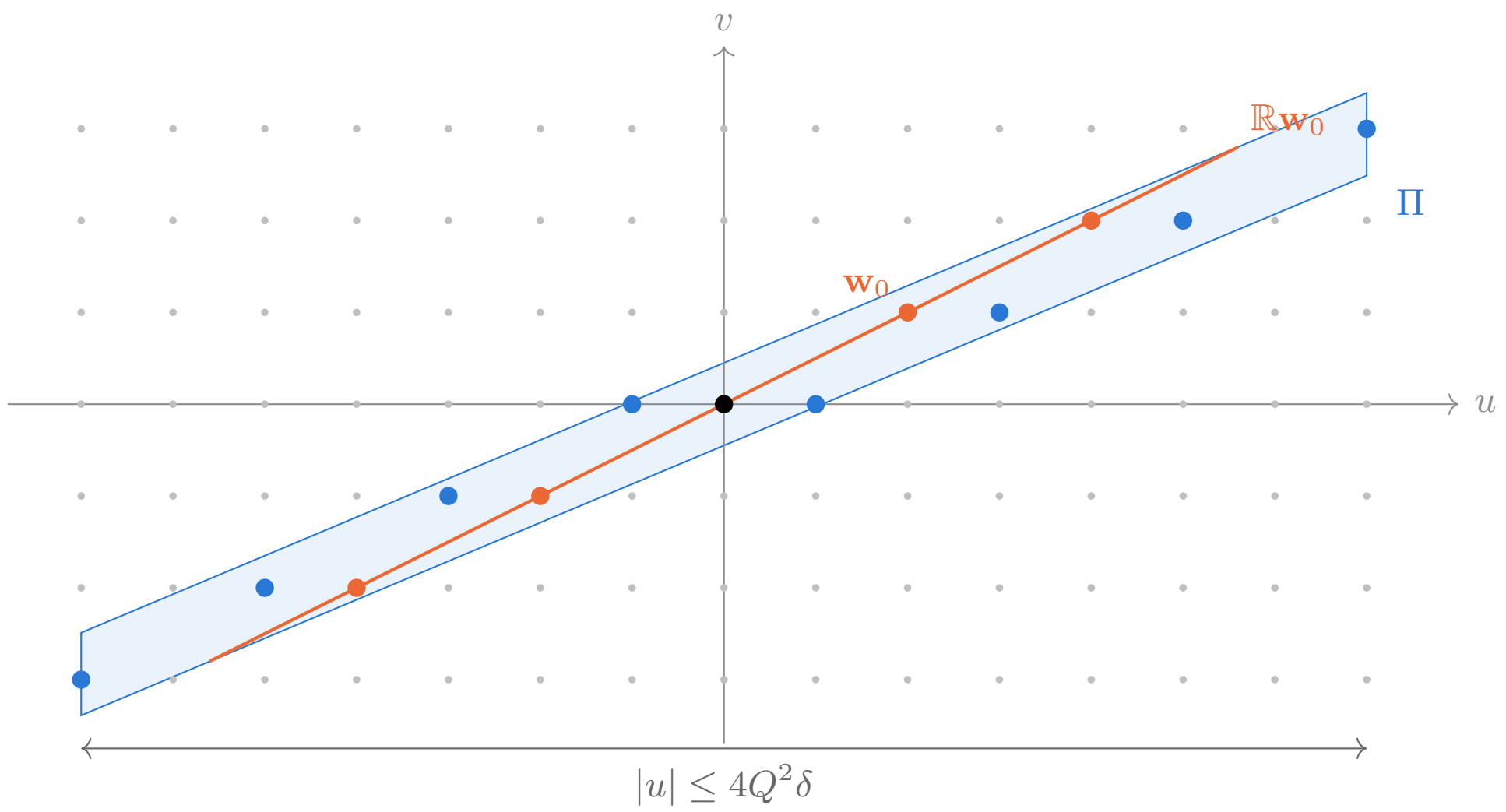

	\smallskip
	\emph{All pairs lie in one parallelogram (Sections~\ref{sec:lattice} and~\ref{sec:eq}).} For every solution, $|u|\le4Q^2\delta$ and $u\alpha-v$ is small. So $(u,v)$ is an integer point of a parallelogram $\Pi$ which depends only on $\alpha$ and the parameters, and its area is $A=64(\kappa\delta^2+T^2\delta)$ (Lemma~\ref{lem:par}). The number of integer points of $\Pi$ is $O(A)$, apart from the points on one line through the origin (Figure~\ref{fig:Pi}). On this line the greatest common divisors $\gcd(u,v)$ are large, and it has to be treated separately: a fraction which has many solutions on the line turns out to have a small value of $|r_0a-b_0q^2|$ for one fixed pair $(r_0,b_0)$, and such fractions are again few by the count of Baier and Zhao.
	
	\smallskip
	\emph{The balance.} Each minor fraction which is not of the exceptional kind has at least $\frac14\delta T$ solutions off the line, and the number of all such solutions is at most $Q^{o(1)}A(1+\mu/T)$. Hence the number of these fractions is at most
	\[
	Q^{o(1)}\,\frac A{\delta T}\Bigl(1+\frac\mu T\Bigr)\ll Q^{o(1)}\Bigl(\frac{\kappa\delta}T+T\Bigr)\Bigl(1+\frac\mu T\Bigr).
	\]
	With $\delta\approx T^{-1/2}$ the two terms in the first bracket are balanced for $T=\kappa^{2/5}$, and the bound becomes $Q^{o(1)}(\kappa^{2/5}+\mu)$. Table~\ref{tab:param} lists the parameters.
	
	\begin{table}[htb]
		\centering\small
		\renewcommand{\arraystretch}{1.25}
		\begin{tabular}{llll}
			\toprule
			parameter & definition & at $N=Q^3$ & role\\
			\midrule
			$\mu$ & $Q^3/N$ & $1$ & average number of fractions in the interval\\
			$\kappa$ & $Q^4/N$ & $Q$ & length of the interval divided by the spacing\\
			$T$ & $\kappa^{2/5}$ & $Q^{2/5}$ & the squares in the equation are $t^2$, $T<t\le2T$\\
			$\delta$ & $T^{-1/2+\eta}$ & $Q^{-1/5+2\eta/5}$ & $|u|\le\delta q^2$; each minor fraction has $\ge\frac12\delta T$ solutions\\
			$D_0$ & $\lfloor T^{5/4}\rfloor$ & $\approx Q^{1/2}$ & Dirichlet parameter for $\theta=\bar a/q^2$\\
			$A$ & $64(\kappa\delta^2+T^2\delta)$ & $Q^{3/5+O(\eta)}$ & area of the parallelogram $\Pi$\\
			\bottomrule
		\end{tabular}
		\caption{The parameters in the proof of Proposition~\ref{thm:thin}.}\label{tab:param}
	\end{table}
	
	\subsection{Quadratic Weyl sums}\label{sec:weyl}
	The aim of this subsection is Lemma~\ref{lem:many}, which says that $\|t^2\theta\|<\delta$ for at least $\frac12\delta T$ integers $t\in(T,2T]$, if $\theta$ has a rational approximation with a denominator between $T$ and $T^{5/4}$. For $\beta\in\R$, $X\ge2$ and $h\in\Z$ we put
	\[
	S_h(\beta;X)=\sum_{X<t\le2X}\e(h\beta t^2).
	\]
	(In this subsection $X$ is a real number, not the indeterminate.) We begin with a standard estimate.
	
	\begin{lemma}\label{lem:summin}
		Let $\gcd(c,s)=1$, $s\ge1$, $|\beta-c/s|\le s^{-2}$, and $J,X\ge1$. Then
		\[
		\sum_{m=1}^J\min\bigl(X,\|m\beta\|^{-1}\bigr)\ll\Bigl(1+\frac Js\Bigr)\bigl(X+s\log2s\bigr).
		\]
	\end{lemma}
	
	\begin{proof}
		For $s\le3$ the bound $JX$ suffices. Otherwise we split $[1,J]$ into at most $1+4J/s$ blocks of at most $\lfloor s/2\rfloor$ consecutive integers. For $m\ne m'$ in one block, $i=|m-m'|$ satisfies $0<i<s/2$, so $\|ic/s\|\ge1/s$ and
		\[
		\|i\beta\|\ge\frac1s-\frac i{s^2}>\frac1{2s}.
		\]
		Thus the points $m\beta$ of a block are $1/(2s)$-separated modulo $1$. At most $3$ of them satisfy $\|m\beta\|<1/(2s)$, and on each side of $0$ the $i$-th of the others has distance at least $i/(2s)$ from $0$. So each block contributes at most $3X+2\sum_{i\le s}2s/i\ll X+s\log2s$.
	\end{proof}
	
	\begin{lemma}[an average of quadratic Weyl sums]\label{lem:weyl}
		Let $X\ge2$, $H\ge1$, $\gcd(c,s)=1$, $1\le s\le X^4$ and $|\beta-c/s|\le s^{-2}$. For every $\rho>0$,
		\[
		\frac1H\sum_{h=1}^H|S_h(\beta;X)|\ \ll_\rho\ (XH)^\rho\Bigl(\sqrt X+\frac X{\sqrt s}+\sqrt{\frac sH}\Bigr).
		\]
	\end{lemma}
	
	\begin{proof}
		By Cauchy's inequality it suffices to bound $\frac1H\sum_h|S_h|^2$ by the square of the right-hand side.
		
		\emph{Weyl differencing.} We expand the square and group the pairs $t<t'$ by $i=t'-t$:
		\[
		|S_h|^2\le X+1+2\sum_{1\le i<X}\Bigl|\sum_t\e\bigl(h\beta(2it+i^2)\bigr)\Bigr|\le X+1+2\sum_{1\le i<X}\min\bigl(X+1,\|2hi\beta\|^{-1}\bigr).
		\]
		Here the inner sum runs over an interval of length at most $X$, and we have used that a geometric sum with ratio $\e(\gamma)$ has modulus at most $\min(\text{length}+1,\|\gamma\|^{-1})$.
		
		\emph{The sum over $h$.} Put $m=2hi\le2HX$. Each $m$ arises from at most $\tau(m)\ll_\rho(XH)^\rho$ pairs $(h,i)$. By Lemma~\ref{lem:summin},
		\[
		\frac1H\sum_{h\le H}|S_h|^2\ll_\rho X+\frac{(XH)^\rho}H\Bigl(1+\frac{HX}s\Bigr)\bigl(X+s\log2s\bigr)\ll_\rho(XH)^{2\rho}\Bigl(X+\frac{X^2}s+\frac sH\Bigr),
		\]
		where we used $\log2s\ll_\rho X^\rho$, which holds as $s\le X^4$, and $X/H\le X$.
	\end{proof}
	
	To detect the condition $\|x\|<\delta$ we use the Fej\'er kernel. For $0<\delta\le\frac12$ let
	\[
	w_\delta(x)=\max\Bigl(0,\,1-\frac{\|x\|}\delta\Bigr).
	\]
	Then $w_\delta$ is at most the indicator function of the set of the $x$ with $\|x\|<\delta$. Moreover $w_\delta$ is $\delta^{-1}$ times the convolution of the indicator function of $[-\delta/2,\delta/2]+\Z$ with itself. Hence $w_\delta(x)=\sum_{h\in\Z}\widehat w_\delta(h)\e(hx)$ with
	\begin{equation}\label{eq:triangle}
		\widehat w_\delta(0)=\delta\qquad\text{and}\qquad0\le\widehat w_\delta(h)=\delta\Bigl(\frac{\sin\pi h\delta}{\pi h\delta}\Bigr)^2\le\min\Bigl(\delta,\frac1{\pi^2\delta h^2}\Bigr)\quad(h\ne0).
	\end{equation}
	
	\begin{lemma}[many small values of $t^2\theta$]\label{lem:many}
		Let $0<\eta<\frac18$. There is $T_0(\eta)$ such that the following holds for $T\ge T_0(\eta)$ and $\delta=T^{-1/2+\eta}$. If $\gcd(c,s)=1$, $T<s\le T^{5/4}$ and $|\theta-c/s|\le s^{-2}$, then
		\[
		\#\bigl\{t\in\Z:\ T<t\le2T,\ \|t^2\theta\|<\delta\bigr\}\ \ge\ \tfrac12\,\delta T .
		\]
	\end{lemma}
	
	\begin{proof}
		Write $S_h=S_h(\theta;T)$. By \eqref{eq:triangle} and $S_{-h}=\overline{S_h}$, the number in question is at least
		\[
		\sum_{T<t\le2T}w_\delta(t^2\theta)=\delta\,\bigl(T+O(1)\bigr)+O\Bigl(\sum_{h\ge1}\widehat w_\delta(h)\,|S_h|\Bigr).
		\]
		The first term is the expected number, and we have to show that the sum over $h$ is smaller.
		
		Let $H_0=\lceil1/\delta\rceil$, so that $1\le\delta H_0\le2$, and fix $\rho=\eta/4$. We split the range of $h$ into $h\le H_0$ and the dyadic ranges $2^{i-1}H_0<h\le2^iH_0$ with $i\ge1$. For $h\le H_0$ we use $\widehat w_\delta(h)\le\delta$, and in the $i$-th dyadic range we use $\widehat w_\delta(h)\le4^{1-i}/(\pi^2\delta H_0^2)$. Lemma~\ref{lem:weyl} with $X=T$ applies, since $s\le T^4$. In its bound, $T/\sqrt s<\sqrt T$ because $s>T$, and
		\[
		\frac s{H_0}\le s\delta\le T^{3/4+\eta}\le T .
		\]
		So the average of $|S_h|$ over $h\le H$ is $\ll(TH)^\rho\sqrt T$ for every $H\ge H_0$, and we obtain
		\[
		\sum_{h\ge1}\widehat w_\delta(h)\,|S_h|\ll\sum_{i\ge0}2^{-i}\,(2^iTH_0)^\rho\sqrt T\ll\Bigl(\frac T\delta\Bigr)^\rho\sqrt T\le T^{1/2+3\rho/2}=T^{1/2+3\eta/8}.
		\]
		This is $o(\delta T)=o(T^{1/2+\eta})$, which proves the lemma for large $T$.
	\end{proof}
	
	\subsection{Integer points in a parallelogram}\label{sec:lattice}
	We need two facts about the integer points of convex sets in the plane. For $z=(x,y)\in\Z^2$ we write $\gcd(z)=\gcd(x,y)$, and $z\ne0$ is called \emph{primitive} if $\gcd(z)=1$.
	
	\begin{lemma}\label{lem:pick}
		\begin{enumerate}[label=\textup{(\alph*)}]
			\item A bounded convex set of area $A$ which contains three non-collinear integer points contains at most $2A+2$ integer points.
			\item If the integer points of an origin-symmetric convex set are collinear, at most two of them are primitive.
		\end{enumerate}
	\end{lemma}
	
	\begin{proof}
		(a) Triangulate the convex hull of the integer points of the set, using all of them as vertices. With $n$ vertices one obtains at least $n-2$ triangles, each of area at least $\frac12$; so $(n-2)/2\le A$.
		
		(b) The integer points lie on a line through $0$, and a line through $0$ contains exactly two primitive integer vectors.
	\end{proof}
	
	\begin{lemma}\label{lem:offline}
		Let $\Pi\subset\R^2$ be a closed parallelogram of area $A$ which is symmetric with respect to the origin and contains a non-zero integer point. Let $\mathbf w_0$ be a non-zero integer point of $\Pi$ which is shortest for the norm whose unit ball is $\Pi$. Then $\mathbf w_0$ is primitive, and
		\[
		\#\bigl((\Pi\setminus\R\mathbf w_0)\cap\Z^2\bigr)\le3A\qquad\text{and}\qquad\sum_{z\in(\Pi\setminus\R\mathbf w_0)\cap\Z^2}\sqrt{\gcd(z)}\le9A .
		\]
	\end{lemma}
	
	\begin{proof}
		The point $\mathbf w_0$ is primitive, because otherwise a proper divisor of it would be a shorter non-zero integer point.
		
		Let $i\ge1$, and suppose that $\Pi$ contains an integer point $iz'$ with $z'\in\Z^2$ and $z'\notin\R\mathbf w_0$. Then $z'\in\Pi/i$, so $\mathbf w_0\in\Pi/i$ by the minimality of $\mathbf w_0$. Thus $\Pi/i$, which has the area $A/i^2$, contains the non-collinear integer points $0$, $\mathbf w_0$ and $z'$. It also contains the convex hull of $\pm\mathbf w_0$ and $\pm z'$, which has area at least $2$; so $A/i^2\ge2$. By Lemma~\ref{lem:pick}(a), the number of such $z'$ is at most
		\[
		\frac{2A}{i^2}+2\le\frac{3A}{i^2}.
		\]
		For $i=1$ this is the first claim. For the second claim we use $\sqrt{\gcd(z)}\le\sum_{i\mid\gcd(z)}\sqrt i$. Each $z$ with $i\mid\gcd(z)$ is of the form $iz'$, so the sum is at most
		\[
		\sum_{i\ge1}\sqrt i\cdot\frac{3A}{i^2}\le9A.\qedhere
		\]
	\end{proof}
	
	\subsection{The count of Baier and Zhao}\label{sec:labels}
	If $b/r$ is a good rational approximation to $\alpha$, the integers $ra-bq^2$ with $(a,q)\in\mathcal S(\alpha)$ are small. The following lemma bounds the number of $(a,q)$ for which $ra-bq^2$ is small, for an arbitrary pair $(r,b)$. It is the counting argument of \cite{BZ}. Recall that $\mu=Q^3/N$.
	
	\begin{lemma}\label{lem:label}
		Let $(r,b)$ be a pair of integers with $\gcd(r,b)=1$, $r\ne0$ and $|r|,|b|\le Q^{O(1)}$, and let $1\le B\le Q^{O(1)}$. Then for every $\alpha$
		\[
		\#\bigl\{(a,q)\in\mathcal S(\alpha):\ |ra-bq^2|\le B\bigr\}\ \le\ Q^{o(1)}\,(1+B+\mu).
		\]
	\end{lemma}
	
	\begin{proof}
		Replacing $(r,b)$ by $(-r,-b)$ we may assume $r>0$. We write $m=ra-bq^2$.
		
		\emph{The pairs with $m=0$.} Then $a/q^2=b/r$, and both fractions are in lowest terms. So $q^2=r$ and $a=b$: there is at most one such pair.
		
		\emph{Localisation of $q$ for a given $m\ne0$.} Since $|a/q^2-\alpha|\le1/N$, the number
		\[
		\frac m{q^2}=r\Bigl(\frac a{q^2}-\alpha\Bigr)+(r\alpha-b)
		\]
		lies within $r/N$ of $r\alpha-b$. On $[Q,2Q]$ the function $x\mapsto m/x^2$ has a derivative of modulus at least $|m|/(4Q^3)$. So $q$ lies in an interval of length at most $8r\mu/|m|$.
		
		\emph{The congruence.} We sort the pairs by $d=\gcd(q,r)$; there are $Q^{o(1)}$ divisors $d$ of $r$. For a given $d$ put
		\[
		g=\gcd(d^2,r),\qquad R=\frac rg,\qquad y=\frac qd .
		\]
		Then $g\mid m$, as $g\mid r$ and $g\mid q^2$, and
		\[
		\frac mg=Ra-b\,\frac{d^2}g\,y^2,\qquad\text{so that}\qquad b\,\frac{d^2}g\,y^2\equiv-\frac mg\pmod R .
		\]
		Here $b$, $d^2/g$ and $y$ are units modulo $R$. Indeed, $\gcd(b,r)=1$; $\gcd(d^2/g,r/g)=1$ by the definition of $g$; and $\gcd(y,r/d)=1$, where $R\mid r/d$ because $d\mid g$. So $y$ is a square root of a given unit modulo $R$, and by Lemma~\ref{lem:roots} it lies in $Q^{o(1)}$ residue classes modulo $R$.
		
		\emph{Counting.} In one of these classes, $q=dy$ runs through an arithmetic progression with the difference $dR$, in an interval of length at most $8r\mu/|m|$. So it takes at most $1+8g\mu/(d|m|)$ values, and $q$ determines $a$. For the given $d$, the admissible $m$ are the multiples of $g$ with $0<|m|\le B$. Their contribution is therefore
		\[
		\le Q^{o(1)}\sum_{\substack{0<|m|\le B\\ g\mid m}}\Bigl(1+\frac{8g\mu}{d|m|}\Bigr)\le Q^{o(1)}\bigl(2B+16\mu\log(3B)\bigr).
		\]
		Summing over the $Q^{o(1)}$ divisors $d$ of $r$ gives the lemma.
	\end{proof}
	
	\subsection{The parameters and the equation}\label{sec:eq}
	From now on we fix $\eps\in(0,\frac18)$, $Q\ge2$, $Q^2\le N\le Q^4$ and $\alpha\in[0,1)$; translating $\alpha$ by an integer changes nothing. We put $\eta=\eps$ and
	\[
	T=\kappa^{2/5},\qquad\delta=T^{-1/2+\eta},\qquad D_0=\lfloor T^{5/4}\rfloor ,
	\]
	as in Table~\ref{tab:param}. By \eqref{eq:Strivial} we may assume that $T\ge T_0$ for any fixed $T_0=T_0(\eps)$. In particular $\delta<\frac14$, $\delta T\ge16$ and $D_0\ge T^{5/4}/2$, and Lemma~\ref{lem:many} applies. As $\kappa\le Q^2$, we have $T\le Q^{4/5}$, so that
	\[
	2T<Q<q\qquad\text{for all }(a,q)\in\mathcal S(\alpha),
	\]
	if $Q\ge T_0$.
	
	\smallskip
	\emph{Major and minor fractions.} For $(a,q)\in\mathcal S(\alpha)$ let $\bar a$ be the inverse of $a$ modulo $q^2$ with $0<\bar a<q^2$, and let
	\[
	\theta=\frac{\bar a}{q^2}.
	\]
	By Dirichlet's approximation theorem, applied to $\theta$, there is a fraction $c/s$ with
	\[
	\gcd(c,s)=1,\qquad1\le s\le D_0,\qquad|s\theta-c|\le\frac1{D_0} .
	\]
	We call $(a,q)$ \emph{major} if $s\le T$ and \emph{minor} if $s>T$.
	
	\smallskip
	\emph{Solutions.} A \emph{solution} is a quintuple $(a,q,t,u,v)$ with $(a,q)\in\mathcal S(\alpha)$ and integers $t,u,v$ such that
	\begin{equation}\label{eq:solution}
		ua-vq^2=t^2,\qquad T<t\le2T,\qquad|u|\le\delta q^2 .
	\end{equation}
	We speak of the solutions \emph{of} $(a,q)$ and of the solutions \emph{of} the pair $(u,v)$.
	
	We first dispose of the major fractions.
	
	\begin{lemma}[major fractions]\label{lem:major}
		The number of major $(a,q)\in\mathcal S(\alpha)$ is at most $Q^{o(1)}(1+T+\mu)$.
	\end{lemma}
	
	\begin{proof}
		The idea is that the approximation $c/s$ to $\bar a/q^2$ produces a pair $(u,v)$ for which $ua-vq^2$ is small, and that only $O(1)$ pairs occur.
		
		\emph{The pair $(u,v)$.} Write $a\bar a=1+\gamma q^2$, and put
		\[
		u=s\bar a-cq^2\qquad\text{and}\qquad v=s\gamma-ca .
		\]
		Then $ua-vq^2=s$. The map $(s,c)\mapsto(u,v)$ is linear with the determinant $-a\bar a+\gamma q^2=-1$, so $\gcd(u,v)=\gcd(s,c)=1$. Also $u\ne0$, since $u\equiv s\bar a\pmod{q^2}$ with $\gcd(\bar a,q)=1$ and $0<s<q^2$.
		
		\emph{Its size.} We have $|u|=q^2|s\theta-c|\le4Q^2/D_0$. Moreover
		\[
		u\alpha-v=u\Bigl(\alpha-\frac a{q^2}\Bigr)+\frac s{q^2},\qquad\text{hence}\qquad|u\alpha-v|\le\frac{4Q^2}{ND_0}+\frac T{Q^2}.
		\]
		So $(u,v)$ lies in an origin-symmetric parallelogram of area
		\[
		4\cdot\frac{4Q^2}{D_0}\Bigl(\frac{4Q^2}{ND_0}+\frac T{Q^2}\Bigr)=\frac{64\kappa}{D_0^2}+\frac{16T}{D_0}\le256+32,
		\]
		where we used $D_0\ge T^{5/4}/2$ and $T^{5/2}=\kappa$. By Lemma~\ref{lem:pick} it contains $O(1)$ primitive integer points.
		
		\emph{Conclusion.} Each major $(a,q)$ satisfies $0<ua-vq^2=s\le T$ for one of these $O(1)$ primitive points $(u,v)$. Their coordinates are at most $4Q^2+1$ in absolute value, so Lemma~\ref{lem:label} with $(r,b)=(u,v)$ and $B=T$ finishes the proof.
	\end{proof}
	
	The next lemma describes the solutions of a fraction.
	
	\begin{lemma}[solutions of a fraction]\label{lem:sol}
		Let $(a,q)\in\mathcal S(\alpha)$ and $T<t\le2T$. There is at most one solution $(a,q,t,u,v)$, and there is one if and only if $\|t^2\theta\|\le\delta$. In every solution $u\ne0$. If $(a,q)$ is minor, it has at least $\frac12\delta T$ solutions.
	\end{lemma}
	
	\begin{proof}
		Let $u_0$ be the least absolute residue of $\bar at^2$ modulo $q^2$, so that $|u_0|=q^2\|t^2\theta\|$, and put $v_0=(u_0a-t^2)/q^2$. This is an integer, because $u_0a\equiv t^2\pmod{q^2}$.
		
		\emph{Uniqueness.} In a solution, $ua\equiv t^2\pmod{q^2}$, hence $u\equiv\bar at^2\equiv u_0\pmod{q^2}$. As $|u|\le\delta q^2<q^2/4$ and $|u_0|\le q^2/2$, this gives $u=u_0$, and then $v=v_0$. Conversely $(a,q,t,u_0,v_0)$ is a solution exactly when $|u_0|\le\delta q^2$, that is, when $\|t^2\theta\|\le\delta$. Moreover $u\ne0$, because $0<t^2<q^2$ and $q^2\mid ua-t^2$.
		
		\emph{Minor fractions.} If $(a,q)$ is minor, then $|\theta-c/s|\le1/(sD_0)\le s^{-2}$ with $T<s\le T^{5/4}$. So Lemma~\ref{lem:many} gives at least $\frac12\delta T$ integers $t\in(T,2T]$ with $\|t^2\theta\|<\delta$.
	\end{proof}
	
	\begin{lemma}[one parallelogram]\label{lem:par}
		For every solution $(a,q,t,u,v)$ the pair $(u,v)$ lies in the origin-symmetric parallelogram
		\[
		\Pi=\Bigl\{(x,y)\in\R^2:\ |x|\le4Q^2\delta,\ \ |x\alpha-y|\le\frac{4Q^2\delta}N+\frac{4T^2}{Q^2}\Bigr\},
		\]
		which has the area
		\[
		A:=64\,\bigl(\kappa\delta^2+T^2\delta\bigr).
		\]
	\end{lemma}
	
	\begin{proof}
		We have $|u|\le\delta q^2\le4Q^2\delta$ and
		\[
		u\alpha-v=u\Bigl(\alpha-\frac a{q^2}\Bigr)+\frac{t^2}{q^2},
		\]
		where $|u(\alpha-a/q^2)|\le4Q^2\delta/N$ and $0<t^2/q^2\le4T^2/Q^2$.
	\end{proof}
	
	\subsection{A pair has few solutions}\label{sec:multiplicity}
	We now bound the number of solutions of a given pair $(u,v)$. The following lemma is the only place, besides Lemma~\ref{lem:label}, where it is used that the moduli are squares. If $\gcd(u,v)=1$ and $\gcd(q,u)=1$, it is the observation from Section~\ref{sec:strategy} that $t/q$ is a square root of $-v$ modulo $|u|$. In general, common factors have to be removed first, and this is the technical part of the proof. Recall that $\mu=Q^3/N$.
	
	\begin{lemma}[square roots in a short interval]\label{lem:mult}
		Let $(u,v)\in\Z^2$ with $u\ne0$ and $|u|,|v|\le Q^{O(1)}$, and let $e_0=\gcd(u,v)$. The number of triples $(a,q,t)$ with $(a,q)\in\mathcal S(\alpha)$, $T<t\le2T$ and $ua-vq^2=t^2$ is at most
		\[
		Q^{o(1)}\Bigl(1+\frac{8\mu\sqrt{e_0}}T\Bigr).
		\]
	\end{lemma}
	
	\begin{proof}
		\emph{Removing common factors.} Write $(u,v)=e_0(r,b)$ with $\gcd(r,b)=1$ and $r\ne0$, so that
		\[
		t^2=e_0\,(ra-bq^2).
		\]
		We fix a divisor $d$ of $|r|$, for which there are $Q^{o(1)}$ choices, and consider the triples with $\gcd(q,|r|)=d$. As in the proof of Lemma~\ref{lem:label} we put
		\[
		g=\gcd(d^2,|r|),\qquad R=\frac{|r|}g,\qquad y=\frac qd .
		\]
		Then $g\mid ra-bq^2$, and $b$, $d^2/g$ and $y$ are units modulo $R$.
		
		\emph{The square.} We have $t^2=e_0g\,m'$ with the positive integer $m'=(ra-bq^2)/g$. Write $e_0g=\sigma\lambda^2$ with $\sigma$ squarefree and $\lambda\ge1$. Then $\lambda\mid t$ and $(t/\lambda)^2=\sigma m'$, so $\sigma\mid t/\lambda$. With the positive integer $j=t/(\sigma\lambda)$ we get
		\[
		m'=\sigma j^2\qquad\text{and}\qquad t=\sigma\lambda j .
		\]
		
		\emph{The congruence.} Reducing $m'$ modulo $R$ gives
		\[
		\sigma\,j^2\equiv-b\,\frac{d^2}g\,y^2\pmod R ,
		\]
		and the right-hand side is a unit. So $\sigma$ and $j$ are units modulo $R$, and
		\[
		(j\bar y)^2\equiv-b\,\frac{d^2}g\,\bar\sigma\pmod R .
		\]
		By Lemma~\ref{lem:roots}, the residue $j\bar y$ modulo $R$ lies in $Q^{o(1)}$ classes. We now count the triples for which $j\bar y$ lies in a given class.
		
		\emph{The fraction $j/y$ is reduced.} Let $p$ be a prime dividing $y$. Then $v_p(q)>v_p(d)=v_p(r)$, and $p\nmid a$. Hence $v_p(ra-bq^2)=v_p(r)$, while $v_p(g)=\min(2v_p(d),v_p(r))=v_p(r)$. So $p\nmid m'$, and therefore $p\nmid j$.
		
		\emph{Spacing.} Two distinct reduced fractions $j_1/y_1$ and $j_2/y_2$ in the same class satisfy $R\mid j_1y_2-j_2y_1\ne0$. So
		\[
		\Bigl|\frac{j_1}{y_1}-\frac{j_2}{y_2}\Bigr|\ge\frac R{y_1y_2}\ge\frac{Rd^2}{4Q^2}.
		\]
		
		\emph{Localisation.} We have
		\[
		\sigma\Bigl(\frac jy\Bigr)^2=\frac{m'}{y^2}=\frac{d^2}g\Bigl(\frac{ra}{q^2}-b\Bigr),
		\]
		and $ra/q^2-b$ lies within $|r|/N$ of $r\alpha-b$. So $(j/y)^2$ lies in a fixed interval of length $2d^2R/(\sigma N)$. Since $t>T$ and $y\le2Q/d$, we have $j/y\ge Td/(2Q\sigma\lambda)$. Now two numbers $x_2>x_1\ge x_0>0$ with $x_2^2-x_1^2\le\Lambda$ satisfy $x_2-x_1\le\Lambda/(2x_0)$. Hence $j/y$ lies in an interval of length at most
		\[
		\frac{2d^2R}{\sigma N}\cdot\frac{Q\sigma\lambda}{Td}=\frac{2dRQ\sqrt{e_0g}}{NT\sqrt\sigma}.
		\]
		
		\emph{Counting.} Dividing the length of this interval by the spacing, we see that each class contains at most
		\[
		1+\frac{8\mu\sqrt{e_0g}}{Td\sqrt\sigma}\ \le\ 1+\frac{8\mu\sqrt{e_0}}T
		\]
		of the fractions $j/y$, as $g\le d^2$ and $\sigma\ge1$.
		
		\emph{Recovery.} The reduced fraction $j/y$ determines $j$ and $y$, hence $q=dy$ and $t=\sigma\lambda j$, and then $a$ from $ua-vq^2=t^2$, as $u\ne0$.
	\end{proof}
	
	\subsection{Proof of Proposition~\ref{thm:thin}}\label{sec:end}
	We keep the notation of Section~\ref{sec:eq}; in particular $A=64(\kappa\delta^2+T^2\delta)$ is the area of the parallelogram $\Pi$ of Lemma~\ref{lem:par}. Let $\Sigma$ be the set of the solutions $(a,q,t,u,v)$ with $(a,q)$ minor. By Lemmas~\ref{lem:sol} and~\ref{lem:par}, every minor $(a,q)$ has at least $\frac12\delta T$ solutions, and the pair $(u,v)$ of every solution is an integer point of $\Pi$ with $u\ne0$. If $\Pi$ contains no integer point other than $0$, then there are no minor $(a,q)$, and Lemma~\ref{lem:major} finishes the proof. Otherwise let $\mathbf w_0=(r_0,b_0)$ be a shortest non-zero integer point of $\Pi$, as in Lemma~\ref{lem:offline}.
	
	\smallskip
	\emph{Solutions off the line $\R\mathbf w_0$.} We apply Lemma~\ref{lem:mult} to each integer point $z=(u,v)$ of $\Pi\setminus\R\mathbf w_0$ with $u\ne0$ and sum over $z$ by Lemma~\ref{lem:offline}. The number of the solutions in $\Sigma$ whose pair $(u,v)$ is not on $\R\mathbf w_0$ is at most
	\begin{equation}\label{eq:off}
		Q^{o(1)}\sum_{z\in(\Pi\setminus\R\mathbf w_0)\cap\Z^2}\Bigl(1+\frac{8\mu\sqrt{\gcd(z)}}T\Bigr)\le Q^{o(1)}A\Bigl(1+\frac\mu T\Bigr).
	\end{equation}
	
	\smallskip
	\emph{Solutions on the line.} If $r_0=0$, no solution has its pair on $\R\mathbf w_0$, because $u\ne0$. Otherwise the integer points on the line are the $i\mathbf w_0$ with $i\in\Z$, since $\mathbf w_0$ is primitive. For a solution with $(u,v)=i\mathbf w_0$ the equation \eqref{eq:solution} reads
	\[
	t^2=i\,m\qquad\text{with}\qquad m=r_0a-b_0q^2\ne0 .
	\]
	The integer $m$ depends on $(a,q)$ only. We call a minor $(a,q)$ \emph{linear} if at least $\frac14\delta T$ of its solutions have their pair on $\R\mathbf w_0$. By Lemma~\ref{lem:sol} these solutions have distinct $t$. Since $m\mid t^2$, all these $t$ are divisible by
	\[
	\beta(m)=\prod_{p^j\|m}p^{\lceil j/2\rceil}.
	\]
	So $\frac14\delta T\le T/\beta(m)+1$. As $\delta T\ge16$, this gives $\beta(m)\le8/\delta$ and
	\[
	|m|\le\beta(m)^2\le64\,\delta^{-2}\le64\,T.
	\]
	Since $\mathbf w_0\in\Pi$ has coordinates $\le Q^{O(1)}$, Lemma~\ref{lem:label} with $(r,b)=(r_0,b_0)$ and $B=64T$ shows that there are at most $Q^{o(1)}(1+T+\mu)$ linear $(a,q)$.
	
	\smallskip
	\emph{Conclusion.} Every minor $(a,q)$ that is not linear has at least $\frac14\delta T$ solutions off the line. By \eqref{eq:off}, the number of these $(a,q)$ is at most
	\[
	Q^{o(1)}\,\frac A{\delta T}\Bigl(1+\frac\mu T\Bigr)\le Q^{o(1)}\Bigl(\frac{\kappa\delta}T+T\Bigr)\Bigl(1+\frac\mu T\Bigr)\le Q^{o(1)}\,T^{\eta}\,(T+\mu),
	\]
	because
	\[
	\frac{\kappa\delta}T=\frac{T^{5/2}\,T^{-1/2+\eta}}T=T^{1+\eta}.
	\]
	Adding the major $(a,q)$ (Lemma~\ref{lem:major}) and the linear ones, we get $\#\mathcal S(\alpha)\le Q^{o(1)}T^\eta(1+T+\mu)$. As $T=\kappa^{2/5}\le Q^{4/5}$ and $\eta=\eps$, this is $\ll_\eps Q^{2\eps}(1+\mu+\kappa^{2/5})$, and $\eps$ is arbitrary. \qed
	
	\subsection{Proof of Theorem~\ref{main:squares}}\label{sec:LSsq}
	Let $\eps>0$. In the notation \eqref{eq:Kdef}, Proposition~\ref{thm:thin} is a bound for $K_{X^2}(1/N,Q)$ when $Q^2\le N\le Q^4$. We first extend it to all ranges: for $2\le Y\le Q$ and every integer $N\ge1$,
	\begin{equation}\label{eq:allY}
		K_{X^2}\Bigl(\frac1N,\,Y\Bigr)\ll_\eps Q^\eps\Bigl(1+\frac{Q^3}N+\Bigl(\frac{Q^4}N\Bigr)^{2/5}\Bigr).
	\end{equation}
	There are three cases.
	\begin{itemize}[leftmargin=2em]
		\item If $N<Y^2$, each $q\le2Y$ has at most $2q^2/N+1\le9Y^2/N$ numerators $a$ with $a/q^2$ in a given interval of length $2/N$. As there are at most $Y+1\le\frac32Y$ values of $q$, there are at most $14Y^3/N$ fractions.
		\item If $N>Y^4$, the spacing $(2Y)^{-4}$ of the fractions leaves at most $1+32Y^4/N\le33$ of them.
		\item If $Y^2\le N\le Y^4$, Proposition~\ref{thm:thin} applies with $Y$ in place of $Q$.
	\end{itemize}
	Each of these bounds is at most the right-hand side of \eqref{eq:allY}, because $Y\le Q$.
	
	Now let $N\ge2$; for $N=1$ the bound $\Delta_{X^2}(Q,1)\ll Q^3$ of \eqref{eq:trivial} suffices. Lemma~\ref{lem:dyadic} with $Y_0=2$ and \eqref{eq:allY} give
	\[
	\Delta_{X^2}(Q,N)\ll N+N\log(2Q)\,Q^\eps\Bigl(1+\frac{Q^3}N+\Bigl(\frac{Q^4}N\Bigr)^{2/5}\Bigr)\ll_\eps Q^{2\eps}\bigl(N+Q^3+Q^{8/5}N^{3/5}\bigr),
	\]
	and $\eps$ is arbitrary. \qed
	
	\subsection{Where squares are used, and what limits the exponent}\label{sec:limits}
	We close with two remarks on the proof.
	
	\smallskip
	\emph{Where squares are used.} Lemma~\ref{lem:label} counts square roots $y=q/d$ of a unit modulo $R$. In Lemma~\ref{lem:mult} the corresponding roots are $j\bar y$, for fixed $(u,v)$ and $d$. Both estimates fail for the moduli $q^2+1$. To see this, take $N=Q^3$, $g=\lceil\sqrt Q\rceil$, $q=gh\in(Q,2Q]$ and $a=h^2$. The $\asymp\sqrt Q$ fractions $a/(q^2+1)$ are reduced, and their distances from $1/g^2$ equal $1/(g^2(q^2+1))\le1/N$; these are the fractions \eqref{eq:cluster}. They all have the same value
	\[
	g^2a-(q^2+1)=-1,
	\]
	whereas the analogue of Lemma~\ref{lem:label} with $B=1$ would give $Q^{o(1)}$ fractions. For the analogue of Lemma~\ref{lem:mult}, fix one integer $t_0\in(T,2T]$ and put $(u,v)=-t_0^2(g^2,1)$. Then
	\[
	ua-v(q^2+1)=t_0^2\qquad\text{and}\qquad\gcd(u,v)=t_0^2.
	\]
	So every fraction in the cluster supplies a triple for this one pair. Since $\mu=1$, the analogue of Lemma~\ref{lem:mult} would give only $Q^{o(1)}(1+8t_0/T)=Q^{o(1)}$ triples.
	
	\smallskip
	\emph{What limits the exponent.} Consider a smaller width $\delta=T^{-\sigma}$ with $\frac12<\sigma<1$. If each minor fraction still supplied $\gg\delta T$ solutions, the count off the shortest line would give
	\[
	Q^{o(1)}\bigl(T+\kappa T^{-1-\sigma}\bigr)\Bigl(1+\frac\mu T\Bigr)
	\]
	fractions with $\gg\delta T$ solutions off that line. Its two terms balance at $\kappa=T^{2+\sigma}$. On the line, however, the argument in Section~\ref{sec:end} gives only $|m|\ll T^{2\sigma}$, so a direct application of Lemma~\ref{lem:label} would cost $Q^{o(1)}(T^{2\sigma}+\mu)$.
	
	One can retain the stronger condition $\beta(m)\le V:=8/\delta$ obtained there. Write $|m|=cj^2$ with $c$ squarefree. Then $\beta(m)=cj$, so at most $2\sum_{c\le V}\lfloor V/c\rfloor\ll V\log(2V)$ non-zero values of $m$ occur. For each fixed $d\mid|r_0|$ in the proof of Lemma~\ref{lem:label}, the sum of the terms $g/|m|$ over these $m$ with $g\mid m$ is $\ll\log(2V)$, since $|m|\le V^2$. Summing the bound for each $m$ in that proof now gives $Q^{o(1)}(V+\mu)$ linear fractions, which is sufficient because $V\ll T$.
	
	At this balance the area in Lemma~\ref{lem:major} is $64\kappa/D_0^2+16T/D_0$. The choice $D_0=\lfloor T^{5/4}\rfloor$ makes the first term of order $T^{\sigma-1/2}$. Taking $D_0=\lfloor T^{1+\sigma/2}\rfloor$ keeps the area bounded, but then the lower bound $\gg\delta T$ for minor fractions is needed throughout $T<s\le D_0$. Such an estimate at the width $T^{-\sigma}$ would give the exponent $1/(2+\sigma)$ at $N=Q^3$. Lemma~\ref{lem:many} gives the exponent $\frac25$, using the wider interval $\delta=T^{-1/2+\eta}$ and the range $T<s\le T^{5/4}$.

	\section*{Acknowledgements}
	The first author thanks the Ramakrishna Mission Vivekananda Educational and Research Institute for excellent working conditions. The second author was supported by the University of New South Wales and the Australian Government through UPA and RTP scholarships, respectively. The third author was supported by grants 207010146582502, 207010146022301, and 207010145609901. This support also enabled the third author to host the first and second authors for research visits. The fourth author was supported by FRG Grant PS71536 at the University of New South Wales.
	
	The authors made extensive use of Anthropic's Claude and OpenAI's ChatGPT in developing the results and proofs and in preparing the manuscript. They have checked the output generated to date for accuracy and plan to refine the presentation shortly.

\end{document}